\documentclass[epsf]{article}
\usepackage{tabularx}
\usepackage{latexsym}
 
\newcommand{\bigzerou}{%
\smash{\lower1.7ex\hbox{\bg 0}}}
\newtheorem{theorem}{Theorem}
\newtheorem{prop}{Proposition}
\newtheorem{defi}{Definition}

\newtheorem{conj}{Conjecture}

\newcommand{\ba}{\begin{eqnarray}}
\newcommand{\ea}{\end{eqnarray}}
\newcommand{\no}{\nonumber}

\def\d{{\partial}}
\newcommand{\mapright}[1]{%
\smash{\mathop{%
\hbox to 1.0cm{\rightarrowfill}}\limits^{#1}}}
\newcommand{\mapleft}[1]{%
\smash{\mathop{%
\hbox to 1.3cm{\leftarrowfill}}\limits^{#1}}}

\begin{document}
\title{
{\bf Multi-Point Virtual Structure Constants and Mirror Computation of $CP^2$-model  }}
\author{Masao Jinzenji ${}^{(1)}$, Masahide Shimizu ${}^{(2)}$ \\
\\
\it (1) Department of Mathematics, Graduate School of Science \\
\it Hokkaido University \\
\it  Kita-ku, Sapporo, 060-0810, Japan\\
{\it e-mail address: jin@math.sci.hokudai.ac.jp}\\
\it (2) Institute for the Advancement of Higher Education  \\
\it Hokkaido University \\
\it  Kita-ku, Sapporo, 060-0817, Japan\\
{\it e-mail address: masa.shimi.2875@gmail.com}}
\maketitle
\begin{abstract}
In this paper, we propose a geometrical approach to mirror computation of genus $0$ Gromov-Witten invariants of $CP^2$.
We use multi-point virtual structure constants, which are defined as intersection numbers of a compact moduli 
space of quasi maps from $CP^1$ to $CP^2$ with $2+n$ marked points. We conjecture that some generating functions 
of them produce mirror map and the others are translated into generating functions of Gromov-Witten invariants 
via the mirror map. We generalize this formalism to open string case. In this case, we have to introduce 
infinite number of deformation parameters to obtain results that agree with some known results of open Gromov-Witten 
invariants of $CP^2$. We also apply multi-point virtual structure constants to compute closed and open Gromov-Witten 
invariants of a non-nef hypersurface in projective space. This application simplifies the computational process 
of generalized mirror transformation. 
\end{abstract}
\section{Introduction}
Mirror computation of the genus $0$ Gromov-Witten invariants of $CP^n$ has been investigated by several authors \cite{bara, iri}.
They were motivated by the Landau-Ginzburg potential proposed by Eguchi, Hori and Xiong in \cite{ehx}:
\begin{eqnarray}
x_1+x_2+\cdots+x_n+\frac{e^t}{x_1 x_2 \cdots x_n}.
\end{eqnarray}
In \cite{bara}, Barannikov introduced a pair $(X,f)$, 
\begin{eqnarray}
X=\{x_0 x_1\cdots x_n=1\}\subset{\bf C}^{n+1},\;f:X\rightarrow {\bf C},\;f=x_0+x_1+\cdots+x_n,
\end{eqnarray}
and considered a potential function,
\begin{equation}
\hat{F}(x;t)=f+\sum_{m=0}^{n}t^m(\sum_{i=0}^{n}x_i)^m,
\end{equation} 
and oscillating integrals,
\begin{equation}
\varphi_{k}(t,\hbar)=\int_{\Delta_k}\exp(\frac{\hat{F}(x;t)}{\hbar})\frac{dx_0 dx_1\cdots dx_n}{d(x_0\cdots x_n)},
\end{equation} 
where $\Delta_k$ is some appropriate relative $n$-cycle of $X$. He next introduced functions $\psi_{k}$ 
satisfying the following normalization condition, 
\begin{equation}
\psi_{k}(t,\hbar)=\sum_{m=0}^{n}u_m(t,\hbar)\hbar\frac{\d\varphi_{k}}{\d t^m}(0,\hbar), 
\end{equation}
and having good asymptotic behavior in $\hbar$. Here $u_m(t,\hbar)=\delta_{m0}+\sum_{j=1}^{\infty}\frac{1}{(\hbar)^j}u_m^{(-j)}(t)$. 
The main result of \cite{bara} is given as follows. If one makes a change of parameters $y^m=u_m^{(-1)}(t)$, then Picard-Fuchs 
equation for $\psi_k$'s takes the form:
\begin{equation}
\frac{\d^2\psi_k}{\d y^i \d y^j}=\frac{1}{\hbar}\sum_{m=0}^{n}A_{ij}^{m}(y)\frac{\d \psi_k}{\d y^m},
\end{equation}  
and 
\begin{equation}
A_{ij}^{n-m}(y)=\d_i\d_j\d_m F_{CP^n}(y),
\end{equation}
where $F_{CP^n}(y)$ is the generating function of genus $0$ Gromov-Witten invariants of $CP^n$.
His method is fundamentally the same as use of Birkhoff factorization to construct connection 
matrices of quantum cohomology from Picard Fuchs equation including $\hbar$ \cite{guest, iri}.
In \cite{iri}, Iritani pointed out this fact and reformulated Barannikov's result in terms of extended 
$I$-function of $CP^n$ \cite{irip}.

In this paper, we propose a geometric approach to mirror computation of genus $0$ Gromov-Witten invariants of $CP^2$
by extending our previous results in \cite{mmg}. In \cite{mmg}, we proposed a general conjecture that for a
toric manifold $X$, we can construct a compact moduli space of quasi maps from $CP^1$ to $X$ of degree ${\bf d}$
with two marked points, which we denote by $\widetilde{Mp}_{0,2}(X,{\bf d})$. We also conjectured 
that an intersection number $w({\cal O}_{\alpha}{\cal O}_{\beta})_{0,{\bf d}}$ of $\widetilde{Mp}_{0,2}(X,{\bf d})$,
which is defined as an analogue of the genus $0$ Gromov-Witten invariant 
$\langle{\cal O}_{\alpha}{\cal O}_{\beta}\rangle_{0,{\bf d}}$ of $X$, gives us the data of the B-model of mirror 
computation. Especially, $w({\cal O}_{1}{\cal O}_{\beta})_{0,{\bf d}}$ gives us the information of 
the mirror map by the following correspondence.
\begin{eqnarray}
t^{\alpha}&=&\eta^{\alpha\beta}\bigl(\sum_{\bf d}w({\cal O}_{1}{\cal O}_{\beta})_{0,{\bf d}}e^{\bf d\cdot x}\bigr),\no\\
          &=&x^{\alpha}+ \eta^{\alpha\beta}\bigl(\sum_{\bf d\neq 0}w({\cal O}_{1}{\cal O}_{\beta})_{0,{\bf d}}e^{\bf d\cdot x}\bigr).
\label{igmir}          
\end{eqnarray} 
In the above formula, $\eta^{\alpha\beta}$ is the inverse of the classical intersection matrix $\eta_{\alpha\beta}=\int_{X}\alpha\wedge\beta$ 
and ${\bf d\cdot x}=\displaystyle{\sum_{\gamma}d_{\gamma}x^{\gamma}}$ where $\gamma$ runs through additive generators of $H^{1,1}(X)$. 
The subscript $\alpha$ in (\ref{igmir}) is not restricted to $H^{1,1}(X)$ and can vary at least the range of the sub-ring of $H^{*,*}(X)$ 
multiplicatively generated by $H^{1,1}(X)$. But the deformation parameters $x^{\gamma}$ are restricted to $H^{1,1}(X)$. 
This restriction comes from the fact that the intersection number $w({\cal O}_{\alpha}{\cal O}_{\beta})_{0,{\bf d}}$ is a 
two point correlation function. In order to include deformation parameters that couple to cohomology elements other than $H^{1,1}(X)$,
we have to define $2+n$ point correlation functions. This task forces us to construct a compact moduli space 
of quasi maps with $2+n$ marked points. In constructing $\widetilde{Mp}_{0,2}(N,d):=\widetilde{Mp}_{0,2}(CP^{N-1},d)$ in \cite{mmg},
the number two is special because it is based on geometric invariant theory of ${\bf C}^{\times}$ action on $CP^1$, that keeps $0$ 
and $\infty$ in $CP^1$ fixed.
On the other hand, Alexeev and Guy considered various compactification of the moduli space of (complex structure of) $CP^1$ with marked 
points in \cite{alex}. In their examples, there exists a compactification that corresponds to ${\bf C}^{\times}$ geometric invariant theory.
We denote this moduli space with $2+n$ marked points by $\overline{M}_{0,2|n}$ because the first two marked points $0$ and $\infty$ are 
special. Main difference from the moduli space $\overline{M}_{0,2+n}$, which is compactified by $PSL(2,{\bf C})$ geometric invariant theory, 
is given as follows. In the open stratum of $\overline{M}_{0,2|n}$, the $n$ marked points are distinct from $0$ and $\infty$, but they can coincide with each other 
in ${\bf C}^{\times}=CP^1-\{0,\infty\}$. Boundaries of $\overline{M}_{0,2|n}$ consist of stable curves of chain shape, whose component $CP^1$'s  
are connected at $0$ and $\infty$. The boundary structure of $\widetilde{Mp}_{0,2}(N,d)$ given in \cite{mmg} is the same as the one of $\overline{M}_{0,2|n}$. 
In this paper, we construct $\widetilde{Mp}_{0,2|n}(N,d)$, the compact moduli space of quasi maps from $CP^1$ to $CP^{N-1}$ with $2+n$ marked points,
by combining the construction processes of $\widetilde{Mp}_{0,2}(N,d)$ with the one of  $\overline{M}_{0,2|n}$.  
We can then define a intersection number $w({\cal O}_{h^a}{\cal O}_{h^b}|\prod_{j=1}^{n}{\cal O}_{h^{m_j}})_{0,d}$, which is an analogue 
of the genus $0$ Gromov-Witten invariant $\langle{\cal O}_{n^a}{\cal O}_{h^b}\prod_{j=1}^{n}{\cal O}_{h^{m_j}}\rangle_{0,d}$. Here $h$ 
is the hyperplane class in $H^{*,*}(CP^{N-1})$. This intersection number satisfies the puncture axiom and the divisor axiom of the 
Gromov-Witten invariant with respect to the operator insertions that correspond to the latter $n$ marked points.
Moreover, we can derive a closed formula of $w({\cal O}_{n^a}{\cal O}_{h^b}|\prod_{j=1}^{n}{\cal O}_{h^{m_j}})_{0,d}$ by applying 
the same localization technique as the one used in \cite{mmg}. With this set-up, we can generalize (\ref{igmir}) in the $CP^{N-1}$ case, 
to include deformation parameters $x^j$ coupled to $h^j$ $(j=0,1,2,\cdots,N-1)$. 
\begin{eqnarray}
t^{i}(x^0,\cdots,x^{N-1})&=&
\eta^{ij}\bigl(\sum_{d=0}^{\infty}\sum_{m_l\geq 0}w({\cal O}_{1}{\cal O}_{h^j}|\prod_{l=0}^{N-1}({\cal O}_{h^l})^{m_l})_{0,d}
\prod_{l=0}^{N-1}\frac{(x^l)^{m_l}}{m_l!}\bigr),\no\\
          &=&x^{i}+ \eta^{ij}\bigl(\sum_{d=1}^{\infty}\sum_{m_l\geq 0}w({\cal O}_{1}{\cal O}_{h^j}|\prod_{l=2}^{N-1}({\cal O}_{h^l})^{m_l})_{0,d}e^{dx^1}
\prod_{l=2}^{N-1}\frac{(x^l)^{m_l}}{m_l!}\bigr).
\label{igmir2}          
\end{eqnarray}
Let $w({\cal O}_{h^a}{\cal O}_{h^b}(x^0,x^1,\cdots,x^{N-1}))_{0}$ be the generating function of 
$w({\cal O}_{h^a}{\cal O}_{h^b}|\prod_{l=0}^{N-1}({\cal O}_{h^l})^{m_l})_{0,d}$. We conjecture that substitution of inversion of (\ref{igmir2}) into 
$w({\cal O}_{h^a}{\cal O}_{h^b}(x^0,x^1,\cdots,x^{N-1}))_{0}$ results in $\langle{\cal O}_{h^a}{\cal O}_{h^b}(t^0,t^1,\cdots,t^{N-1})\rangle_{0}$, the 
generating function of $\langle{\cal O}_{h^a}{\cal O}_{h^b}\prod_{l=0}^{N-1}({\cal O}_{h^l})^{m_l}\rangle_{0,d}$.
In this paper, we test numerically this conjecture in the $CP^2$ case. The result indeed supports our conjecture. 
At this stage, a natural question may occur. Does the mirror map constructed in (\ref{igmir2}) coincide with the mirror map obtained from the 
method of Barannikov and Iritani? In this paper,we compare our mirror map for $CP^2$-model with the mirror map derived from Iritani's extended $I$-function.
We find that they do not coincide. But, we cannot conclude only from this fact that our construction has no connection with the standard mirror computation 
of $CP^2$-model. According to Iritani \cite{irip}, there are infinitely many ways to choose a polynomial that couple to the deformation parameter $x^2$.
If we change the polynomial, the mirror map obtained from Birkhoff factorization may vary. Up to now, we have tested only one choice. Therefore, we cannot 
deny the possibility of appropriate choice that reproduces our mirror map.

Next, we apply the moduli space $\widetilde{Mp}_{0,2|2n}(3,2d-1)$ to compute open Gromov-Witten invariants of $CP^2$. We define 
an anti-holomorphic involution of $\widetilde{Mp}_{0,2|2n}(3,2d-1)$ that exchange the two special marked points. We denote 
by $\widetilde{Mp}_{D,1|n}(CP^2/RP^2,2d-1)$ the subset of $\widetilde{Mp}_{0,2|2n}(3,2d-1)$ invariant under the involution.@
It is nothing but the moduli space of quasi maps from disk $\{|z|\leq 1\}$ to $CP^2$. Here, boundary of the disk is 
mapped to $RP^2$, the real Lagrangian submanifold of $CP^2$. The special marked point is $0$ and the remaining $n$ marked points 
can lie freely on $\{0<|z|\leq 1\}$. With this set-up, we compute a intersection number 
$w({\cal O}_{h^a}|\prod_{j=1}^{n}{\cal O}_{h^{m_j}})_{disk,2d-1}$, which is an analogue of the open Gromov-Witten invariant  
$\langle{\cal O}_{h^a}\prod_{j=1}^{n}{\cal O}_{h^{m_j}}\rangle_{disk,2d-1}$. This line of construction is a generalization 
of our previous work \cite{opv} on one point open Gromov-Witten invariants to $1+n$ point Gromov-Witten invariants. 
With this construction, we can compute $w({\cal O}_{h^a}(x^0,x^1,x^2))_{disk}$, the generating function of 
$w({\cal O}_{h^a}|\prod_{l=0}^{2}({\cal O}_{h^l})^{m_l})_{disk,2d-1}$. Then our question is the following.
If we substitute inversion of the mirror map given in (\ref{igmir2}) into $w({\cal O}_{h^a}(x^0,x^1,x^2))_{disk}$, 
can we obtain the generating function $\langle{\cal O}_{h^a}(t^0,t^1,t^2)\rangle_{disk}$? The answer turns out to be "no". 
The result so obtained did not reproduce $\langle{\cal O}_{h^a}\prod_{l=0}^{2}({\cal O}_{h^l})^{m_l}\rangle_{disk,2d-1}$
with lower $d$ that can be computed from localization theorem (naively) applied to the moduli space of stable maps for open 
Gromov-Witten invariants. After some try and errors, we found a way to remedy this disagreement. It is to introduce 
unnatural deformation parameters $t^j$ $(j=3,4,\cdots)$ that couple to $h^j$.  Note that $h^j=0$ $(j\geq 3)$ in $H^{*,*}(CP^2)$, 
but we can formally compute non-zero $w({\cal O}_{h^a}|\prod_{l=0}^{\infty}({\cal O}_{h^l})^{m_l})_{disk,2d-1}$ because
the closed formula for this intersection number is represented in the form of a residue integral. In order to obtain 
the mirror map, we also have to compute $w({\cal O}_{1}{\cal O}_{h^a}|\prod_{l=0}^{\infty}({\cal O}_{h^l})^{m_l})_{0,d}$ with $a\leq -1$
on the closed string side. But it is possible since the closed formula for $w({\cal O}_{h^a}{\cal O}_{h^b}|\prod_{j=1}^{n}{\cal O}_{h^{m_j}})_{0,d}$
is also written in the form of a residue integral. After all, what we have obtained is the following table of open Gromov-Witten invariants 
$\langle({\cal O}_{h^2})^{3d-2}\rangle_{disk,2d-1}$.

\begin{center}
\begin{tabularx}{13cm}{|X|c|c|c|c|c|c|}
\hline\multicolumn{7}{|c|}{\bf Disk Gromov-Witten Invariants of $CP^{2}$ with Maximal $h^{2}$-Insertions}\\
\hline
$d$&$1$&$2$&$3$&$4$&$5$&$6$\\ \hline
$\langle({\cal O}_{h^2})^{3d-2}\rangle_{disk,2d-1}$&$2$&$-\frac{9}{4}$&$\frac{3361}{32}$&
$ -\frac{5784805}{256}$&$\frac{28104787833}{2048}$&$-\frac{291021328876469}{16384}$\\
\hline
\end{tabularx}
\end{center}
The first four numbers coincide with the numerical results obtained from localization computation applied to the moduli space 
of stable maps for open Gromov-Witten invariants.

Our construction of multi-point intersection numbers for quasi maps can be easily generalized to degree $k$ hypersurface 
in $CP^{N-1}$ (we denote it by $M_{N}^{k}$). In this case, we can use our new intersection numbers to 
simplify computational process of generalized mirror transformation \cite{gene, gene0}. In the $k>N$ case, generalized mirror 
transformation for two point Gromov-Witten invariants includes multi-point intersection numbers. 
But in \cite{gene}, the intersection numbers we have as initial data are the ones that correspond to $w({\cal O}_{h^a}{\cal O}_{h^b})_{0,d}$.
Therefore, in order to obtain multi-point intersection numbers, we have to use associativity equation \cite{km}. This 
process made the computation awfully complicated. In the open string case \cite{opv}, this obstacle also made the mirror computation 
of open Gromov-Witten invariants of $M_{N}^{k}$ $(K>N)$ incomplete. Now that we have multi-point intersection numbers on the B-model side,
we can execute the generalized mirror transformation only by one process of coordinate change.   
 In this paper, we apply this idea to compute both closed and open Gromov-Witten invariants of $M_{8}^{9}$. 
It works well as we expected. Especially in the open string case, we don't have to introduce unnatural deformation parameters
in contrast to the $CP^2$ case.  

This paper is organized as follows.
In Section 2, we first construct $\widetilde{Mp}_{0,2|n}(N,d)$, the compact moduli space of quasi maps from $CP^1$ to $CP^{N-1}$ of 
degree $d$ with $2+n$ marked points. After giving definition of the intersection number 
$w({\cal O}_{h^a}{\cal O}_{h^b}|\prod_{j=1}^{n}{\cal O}_{h^{m_j}})_{0,d}$, we derive a closed formula for it by using 
localization technique. Next, we present numerical results of these intersection numbers for $CP^2$ and test our conjecture 
on the mirror computation of $CP^2$-model. In the last part of this section, we compare our mirror map with the one obtained 
from Iritani's extended $I$-function. In Section 3, we construct the moduli space $\widetilde{Mp}_{D,1|n}(CP^2/RP^2,2d-1)$ 
and derive a closed formula for the intersection number $w({\cal O}_{h^a}|\prod_{j=1}^{n}{\cal O}_{h^{m_j}})_{disk,2d-1}$
along the same line as Section 2. Next, we exhibit explicit numerical results of the mirror computation of open Gromov-Witten 
invariants of $CP^2$-model. In Section 4, we apply the formalism in Section 2 and Section 3 to compute closed and open 
Gromov-Witten invariants of $M^{9}_{8}$. We show that our new formalism simplifies the process of generalized mirror 
transformation in \cite{gene} and \cite{opv}, by presenting explicit numerical data. 

\vspace{1cm}

{\bf Acknowledgment} We would like to thank Prof. Yongbin Ruan and Prof. Satoshi Minabe for valuable discussions.
We also thank Prof. Hiroshi Iritani for explaining us his way of mirror computation of the $CP^2$ model.
The research of M.J. is partially supported 
by JSPS grant No. 25400061. 

\section{Closed String Case}

\subsection{Construction of the Moduli Space $\widetilde{Mp}_{0,2|n}(N,d)$}
In this section, we construct $\widetilde{Mp}_{0,2|n}(N,d)$, the 
moduli space of quasi maps from a certain class of semi-stable genus $0$ curves to $CP^{N-1}$ with $2+n$ marked points.
The semi-stable curve we use is a chain of several $CP^1$'s , each of which is connected at $0$ and $\infty$.
Let $l$ be the number of $CP^1$'s in the semi-stable curve. We represent here the stable curve as $\cup_{i=1}^{l}(CP^1)_i$.
Intersection of $(CP^1)_{i}$ and $(CP^1)_{i+1}$ is given by $(\infty)_i=(0)_{i+1}$. $(0)_1$ and $(\infty)_l$ are the 
two special marked points of the semi-stable curve.  A quasi map $\varphi$ from $CP^{1}$ to $CP^{N-1}$ of degree d is defined by,
\begin{equation}
\varphi(s:t)=[\sum_{i=0}^{d}{\bf a}_i s^{i}t^{d-i}],
\label{qm}
\end{equation}
where $s$ and $t$ are homogeneous coordinates $(s:t)$ of $CP^1$.  $[*]$ denotes equivalence class of $*$ under projective equivalence of $CP^{N-1}$.
Our construction is based on construction of $\widetilde{Mp}_{0,2}(N,d)$, the 
moduli space of quasi maps from the above semi-stable curve to $CP^{N-1}$ with two marked points.
Its construction was given in our previous work \cite{mmg}. In constructing $\widetilde{Mp}_{0,2}(N,d)$,
we considered a chain of quasi maps $\cup_{i=1}^{l}\varphi_i$ to represent a quasi map from $\cup_{i=1}^{l}(CP^1)_i$ to $CP^{N-1}$.
The chain  $\cup_{i=1}^{l}\varphi_i$ is classified by ordered partition $(d_1,d_2,\cdots,d_l)$ of $d$ where $d_i$ is 
degree of $\varphi_i$ and it was used to compactify the moduli space.
Then what we have to do in addition to construct $\widetilde{Mp}_{0,2|n}(N,d)$ is to distribute $n$ marked points 
on the semi-stable curve $\cup_{i=1}^{l}(CP^1)_i$. But there occurs one subtlety.
We have to consider the case when some of the $n$ marked points are located at $(0)_i$ or $(\infty)_i$.
To describe this situation, we need to insert $\overline{M}_{0,2|n}$, the moduli space of complex structure 
of $CP^1$ with $2+n$ marked points compactified by $C^{\times}$ geometric invariant theory, to these points. 

Now let us begin construction of $\widetilde{Mp}_{0,2|n}(N,d)$. 
As the first step, we define $Mp_{0,2|n}(N,d)$, which is Zariski-open subset (or bulk part) of $\widetilde{Mp}_{0,2|n}(N,d)$.
It is defined as follows:
\begin{eqnarray}
&&Mp_{0,2|n}(N,d):=\{\bigl(({\bf a}_{0},{\bf a}_{1},\cdots,{\bf a}_{d}),(z_1,z_2,\cdots,z_n)\bigr)\;|\;{\bf a}_i\in{\bf C}^{N},\;{\bf a}_0,{\bf a}_d\neq{\bf 0},\;z_j\in{\bf C}^{\times}\;\}/({\bf C}^{\times})^2.
\label{mp2n}
\end{eqnarray}
In this definition, the first two marked points that correspond to $2$ in the subscript $2|n$ are $0$ and $\infty$ in $CP^1$. We note here that 
$z_j$'s need not be distinct points in ${\bf C}^{\times}=CP^1\setminus\{0,\infty\}$ when $n>1$. We define here the $({\bf C}^{\times})^2$ action in (\ref{mp2n}). Let $(\lambda,\mu)\in ({\bf C}^{\times})^2$.  Then it is given by,
\begin{eqnarray}
&&\lambda\cdot\bigl(({\bf a}_{0},{\bf a}_{1},\cdots,{\bf a}_{d}),(z_1,z_2,\cdots,z_n)\bigr)
=((\lambda{\bf a}_{0},\lambda{\bf a}_{1},\cdots,\lambda{\bf a}_{d}),(z_1,z_2,\cdots,z_n)\bigr),\no\\
&&\mu\cdot\bigl(({\bf a}_{0},{\bf a}_{1},\cdots,{\bf a}_{d}),(z_1,z_2,\cdots,z_n)\bigr)
=(({\bf a}_{0},\mu{\bf a}_{1},\cdots,\mu^{i}{\bf a}_i,\cdots\mu^{d}{\bf a}_{d}),(\mu^{-1} z_1,\mu^{-1} z_2,\cdots,\mu^{-1} z_n)\bigr).
\label{action}
\end{eqnarray}
The first two marked points are $(0:1)$ and $(1:0)$ and the remaining $n$ marked 
points are given by $z_j=\frac{s_j}{t_j}\;(j=1,2,\cdots,n)$ respectively. The condition ${\bf a}_0, {\bf a}_d\neq {\bf 0}$ 
guarantees that images 
$\varphi(0:1)$ and $\varphi(1:0)$ are well-defined in $CP^{N-1}$. 
The first ${\bf C}^{\times}$ action in (\ref{action}) corresponds to 
 projective equivalence of $CP^{N-1}$ and the second one corresponds to automorphism group of $CP^1$ that 
fixes the first two marked points. 
 
When we consider the moduli space of quasi maps with marked points, it is important 
to consider evaluation maps. For the first two marked points, we define evaluation maps $ev_{0}$ and $ev_{\infty}$ as follows:
\begin{eqnarray}
&&ev_{0}([\bigl(({\bf a}_{0},{\bf a}_{1},\cdots,{\bf a}_{d}),(z_1,z_2,\cdots,z_n)\bigr)]):=[{\bf a}_{0}],\no\\
&&ev_{\infty}([\bigl(({\bf a}_{0},{\bf a}_{1},\cdots,{\bf a}_{d}),(z_1,z_2,\cdots,z_n)\bigr)]):=[{\bf a}_{d}],
\label{ev1}
\end{eqnarray}
where $[*]$'s in the l.h.s.'s denote the equivalence classes under the $({\bf C}^{\times})^2$ action. To define evaluation 
maps $ev_{i}\;(i=1,\cdots,n)$ that come from evaluation of a quasi map at $z_i$, 
we have to take care of a subtlety arising from "freckled instantons". A freckled instanton 
is a quasi map $\varphi(s:t)$ whose defining vector valued polynomial $\sum_{i=0}^{d}{\bf a}_i s^{i}t^{d-i}$ is factored as 
follows:
\begin{eqnarray}
\sum_{i=0}^{d}{\bf a}_i s^{i}t^{d-i}=\bigl(\prod_{l=1}^{m}(s-\alpha_l t)\bigr)\cdot(\sum_{j=0}^{d-m}{\bf b}_{j}s^{i}t^{d-m-j}),\;\;(m\geq 1,\;\;
\alpha_l\in{\bf C}^{\times}).
\label{fr}
\end{eqnarray}
In (\ref{fr}), $\sum_{j=0}^{d-m}{\bf b}_{j}s^{j}t^{d-m-j}$ is an irreducible vector valued polynomial. We cannot 
define image of $\varphi$ at $(\alpha_l:1)$ because $\sum_{i=0}^{d}{\bf a}_i(\alpha_l)^{i}={\bf 0}$. But images of 
the other points in $CP^{1}$ are given by $[\sum_{j=0}^{d-m}{\bf b}_{j}s^{j}t^{d-m-j}]$. Therefore, we define $ev_{i}$ 
as follows:
\begin{eqnarray}
&&ev_{i}([\bigl(({\bf a}_{0},{\bf a}_{1},\cdots,{\bf a}_{d}),(z_1,z_2,\cdots,z_n)\bigr)]):=[\sum_{j=0}^{d-m}{\bf b}_{j}(z_i)^{j}].
\label{ev2}
\end{eqnarray} 
If $\sum_{i=0}^{d}{\bf a}_i s^{i}t^{d-i}$ is irreducible, we define,
\begin{eqnarray}
&&ev_{i}([\bigl(({\bf a}_{0},{\bf a}_{1},\cdots,{\bf a}_{d}),(z_1,z_2,\cdots,z_n)\bigr)]):=[\sum_{j=0}^{d}{\bf a}_{j}(z_i)^{j}].
\label{ev3}
\end{eqnarray}

As was suggested in the construction of $\widetilde{Mp}_{0,2}(N,d)$ in \cite{mmg}, we can easily see that $Mp_{0,2|n}(N,d)$ is not a compact space. Therefore, 
we compactify it by adding boundary strata. We start this process by taking $M_{0,2|n}(N,1)$ as an example. In this case, 
it is obtained by taking $({\bf C}^{\times})^2$ quotients of the set $\{\bigl(({\bf a}_{0},{\bf a}_1),(z_1,z_2,\cdots,z_n)\bigr)\;|
\;{\bf a}_{0},{\bf a}_1\neq {\bf 0},z_i\in {\bf C}^{\times} \}$. Using the $({\bf C}^{\times})^2$ action, we can see,
\begin{equation}
Mp_{0,2|n}(N,1)=\{\bigl(([{\bf a}_0],[{\bf a}_1]),(z_1,\cdots, z_n) \bigr)\;|\;z_i\in{\bf C}^{\times}\}=CP^{N-1}\times CP^{N-1}\times({\bf C}^{\times})^n.
\label{com1}               
\end{equation}
Therefore, it is not compact. To compactify this space, we have to add boundaries to describe the situation that some $z_i$'s go to 
$0$ or $\infty$ in $CP^1$. For this purpose, we introduce $\overline{M}_{0,2|n}$, the moduli space of complex structure of $CP^1$
with $n$ marked points compactified by stable curves of chain shape. To construct $\overline{M}_{0,2|n}$, we start from 
an open stratum $M_{0,2|n}$ defined by,
\begin{equation}  
M_{0,2|n}:=\{(z_1,z_2,\cdots,z_n)\;|\;z_i\in {\bf C}^{\times}(=CP^1\setminus \{0,\infty\})\;\}/{\bf C}^{\times},
\end{equation}
where ${\bf C}^{\times}$ action is given by $\mu\cdot(z_1,z_2,\cdots,z_n):=(\mu z_1,\mu z_2,\cdots,\mu z_n)$.
The two marked points that correspond to $2$ in the subscript $2|n$ are $0$ and $\infty$. 

We can easily see that $M_{0,2|1}$ is just an one point set. Then we decompose the subscript set $\{1,2,\cdots,n\}$
into disjoint union of ordered $l$ subsets $(1\leq l\leq n)$:
\begin{equation}
\coprod_{j=1}^{l}A_{j}=\{1,2,\cdots,n \},\;A_i\neq\emptyset.
\end{equation} 
With this set-up, $\overline{M}_{0,2|n}$ is give by a disjoint union of strata as follows:
\begin{equation}
\overline{M}_{0,2|n}=\coprod_{\coprod_{j=1}^{l}A_{j}=\{1,2,\cdots,n \}}\biggl(M_{0,2||A_1|}\times M_{0,2||A_2|}\times\cdots\times M_{0,2||A_l|}\biggr).
\label{m02n}
\end{equation}
In (\ref{m02n}), a point in the stratum labeled by $\coprod_{j=1}^{l}A_{j}$ corresponds to the stable curve $\cup_{i=1}^{l}(CP^1)_i$
where $(CP^1)_i$ has $|A_i|$ marked points $(z_j)_i\;(j\in A_i)$ other than $0$ and $\infty$. 
With this stable curve in mind, how to glue these strata is obvious. Next, we introduce a compact space $\widetilde{M}_{0,2|n}(N,0)$ given by, 
\begin{equation}
\widetilde{Mp}_{0,2|n}(N,0):=CP^{N-1}\times \overline{M}_{0,2|n}.
\end{equation}
For brevity, we also  introduce an auxiliary space $\widetilde{Mp}_{0,2|0}(N,0):=CP^{N-1}$. We then turn back to compactification of $Mp_{0,2|n}(N,1)$.
We decompose the subscript set $\{1,2,\cdots,n\}$ into disjoint union of three subsets:
\begin{equation}
A_{0}\coprod B_1\coprod A_{1}=\{1,2,\cdots,n\}.
\label{dec}
\end{equation}
In (\ref{dec}), each subset can be an empty set. Then $\widetilde{Mp}_{0,2|n}(N,1)$, the compactification of $Mp_{0,2|n}(N,1)$
is given as follows:
\begin{equation}
\widetilde{Mp}_{0,2|n}(N,1)=\coprod_{A_{0}\coprod B_1\coprod A_{1}=\{1,2,\cdots,n \}}
\biggl(\widetilde{Mp}_{0,2||A_{0}|}(N,0)\mathop{\times}_{CP^{N-1}}Mp_{0,2||B_1|}(N,1)\mathop{\times}_{CP^{N-1}}
\widetilde{Mp}_{0,2||A_{1}|}(N,0)\biggr)
\label{wmp1}
\end{equation}
In (\ref{wmp1}), the stratum labeled by $A_{0}\coprod B_1\coprod A_{1}$ corresponds to the configuration 
of marked points where $z_j\;(j\in A_{0})$ (resp. $z_j\;(j\in A_{1})$) goes to $0$ (resp. $\infty$). 
$\widetilde{Mp}_{0,2||A_{0}|}(N,0)\mathop{\times}_{CP^{N-1}}Mp_{0,2||B_1|}(N,1)$ is a fiber product 
with respect to the projection $\pi:\widetilde{Mp}_{0,2||A_{0}|}(N,0)\rightarrow CP^{N-1}$ and 
$ev_{0}:Mp_{0,2||B_1|}(N,1)\rightarrow CP^{N-1}$. $Mp_{0,2||B_{1}|}(N,1)\mathop{\times}_{CP^{N-1}}\widetilde{Mp}_{0,2||A_{1}|}(N,0)$
is also a fiber product 
with respect to $ev_{\infty}:Mp_{0,2||B_{1}|}(N,1)\rightarrow CP^{N-1}$ and the projection
$\pi:\widetilde{Mp}_{0,2||A_{1}|}(N,0)\rightarrow CP^{N-1}$.

We now turn into construction of $\widetilde{Mp}_{0,2|n}(N,d)$. For this purpose, we look back at the construction of 
$\widetilde{Mp}_{0,2}(N,d)$, i.e., the $n=0$ case. In this case, we introduce ordered partition of the degree $d$:
\begin{equation}
OP_{d}:=\{(d_1,d_2,\cdots,d_l)\;|\;d_1+d_2+\cdots+d_l=d,\;d_j\geq 1\;,1\leq l\leq d\;\}.
\end{equation}
In \cite{mmg}, $\widetilde{Mp}_{0,2}(N,d)$ was constructed as follows:
\begin{equation}
\widetilde{Mp}_{0,2}(N,d):=\coprod_{(d_1,\cdots,d_l)\in OP_d}\biggl(Mp_{0,2}(N,d_1)\mathop{\times}_{CP^{N-1}}Mp_{0,2}(N,d_2)\mathop{\times}_{CP^{N-1}}
\cdots\mathop{\times}_{CP^{N-1}}Mp_{0,2}(N,d_l)\biggr),
\label{mpd}
\end{equation}
where $Mp_{0,2}(N,d)$ is the space $Mp_{0,2|0}(N,d)$ defined in (\ref{mp2n}). In (\ref{mpd}), $Mp_{0,2}(N,d_i)\mathop{\times}_{CP^{N-1}}Mp_{0,2}(N,d_{i+1})\;(1\leq i\leq l-1)$
is a fiber product with respect to $ev_{\infty}:Mp_{0,2}(N,d_i)\rightarrow CP^{N-1}$ and $ev_{0}:Mp_{0,2}(N,d_{i+1})\rightarrow CP^{N-1}$. 
A point in the stratum labeled by $(d_1,\cdots,d_l)$ is a chain of quasi maps $\cup_{i=1}^{l}\varphi_i$ where degree of $\varphi_i$ is given by $d_i$.
Fiber products are used to guarantee $\varphi((\infty)_i)=\varphi_{i+1}((0)_{i+1})$.
Construction of 
$\widetilde{Mp}_{0,2|n}(N,d)$ is done by combining (\ref{wmp1}) and (\ref{mpd}). For an ordered partition $(d_1,d_2,\cdots,d_l)\in OP_{d}$, we consider 
ordered decomposition of the subscript set $\{1,2,\cdots,n\}$:
\begin{equation}
\bigl(\coprod_{i=0}^{l}A_i\bigr)\coprod\bigl(\coprod_{i=1}^{l}B_i\bigr)=\{1,2,\cdots,n\},
\label{decl}
\end{equation}
where $A_i$ and $B_j$ can be empty sets. Then $\widetilde{Mp}_{0,2|n}(N,d)$ is given as follows:
\begin{eqnarray}
&&\widetilde{Mp}_{0,2|n}(N,d):=
\coprod_{(d_1,\cdots,d_l)\in OP_d}\coprod_{\bigl(\coprod_{i=0}^{l}A_i\bigr)\coprod\bigl(\coprod_{i=1}^{l}B_i\bigr)=\{1,2,\cdots,n\}}
\biggl(\widetilde{Mp}_{0,2||A_{0}|}(N,0)\mathop{\times}_{CP^{N-1}}Mp_{0,2||B_1|}(N,d_1)\mathop{\times}_{CP^{N-1}}\no\\
&&\mathop{\times}_{CP^{N-1}}\widetilde{Mp}_{0,2||A_{1}|}(N,0)\mathop{\times}_{CP^{N-1}}
Mp_{0,2||B_2|}(N,d_2)\mathop{\times}_{CP^{N-1}}\widetilde{Mp}_{0,2||A_{2}|}(N,0)\mathop{\times}_{CP^{N-1}}
Mp_{0,2||B_3|}(N,d_3)\mathop{\times}_{CP^{N-1}}
\cdots\no\\
&&
\cdots
\mathop{\times}_{CP^{N-1}}\widetilde{Mp}_{0,2||A_{l-1}|}(N,0)
\mathop{\times}_{CP^{N-1}}Mp_{0,2||B_l|}(N,d_l)\mathop{\times}_{CP^{N-1}}\widetilde{Mp}_{0,2||A_{l}|}(N,0)\biggr).
\label{mpdn}
\end{eqnarray}
In (\ref{mpdn}), $Mp_{0,2||B_i|}(N,d_i)\mathop{\times}_{CP^{N-1}}\widetilde{Mp}_{0,2||A_{i}|}(N,0)\mathop{\times}_{CP^{N-1}}
Mp_{0,2||B_{i+1}|}(N,d_{i+1})$ means just \\
$$Mp_{0,2||B_i|}(N,d_i)\mathop{\times}_{CP^{N-1}}
Mp_{0,2||B_{i+1}|}(N,d_{i+1})$$ if $A_i=\emptyset$. Otherwise, it means $$\biggl((Mp_{0,2||B_i|}(N,d_i)\mathop{\times}_{CP^{N-1}}
Mp_{0,2||B_{i+1}|}(N,d_{i+1})\biggr)\times \overline{M}_{0,2||A_{i}|}.$$
A point in the stratum in (\ref{mpdn}) labeled by $(d_1,\cdots,d_l)$ represents a chain of quasi maps 
$\varphi_{i}:(CP^{1})_i\rightarrow CP^{N-1}\;(i=1,2,\cdots,l)$ that satisfies $\varphi_{i}((\infty)_i)=\varphi_{i+1}((0)_{i+1})\;(1\leq i\leq l-1)$.
If it is in the stratum labeled by $\bigl(\coprod_{i=0}^{l}A_i\bigr)\coprod\bigl(\coprod_{i=1}^{l}B_i\bigr)=\{1,2,\cdots,n\}$,
the marked point $z_{j}\;(j\in B_i)$ lies inside ${\bf C}^{\times}\subset (CP^1)_{i}$ and the marked point $z_{j}\;(j\in A_i)$
are mapped to $\varphi_{i}((\infty)_i)=\varphi_{i+1}((0)_{i+1})\;\;(i=1,\cdots,l-1)$ (resp. $\varphi_1((0)_1)$ if $i=0$ and 
$\varphi_l((\infty)_l)$ if $i=l$). With this set-up, we can easily extend the definition of the evaluation maps $ev_{0}, ev_{\infty},
ev_i\;(i=1,2,\cdots,n)$ to whole $\widetilde{Mp}_{0,2|n}(N,d)$.

\subsection{Localization Computation}
In this section, we compute an intersection number $w({\cal O}_{h^a}{\cal O}_{h^b}|\prod_{j=1}^{n}{\cal O}_{h^{m_j}})_{0,d}$ 
on $\widetilde{Mp}_{0,2|n}(N,d)$ by using localization technique developed in \cite{mmg}.
Here, $h$ is the hyperplane class in $H^{*}(CP^{N-1})$.
It is defined by the following formula:
\begin{equation}
w({\cal O}_{h^a}{\cal O}_{h^b}|\prod_{j=1}^{n}{\cal O}_{h^{m_j}})_{0,d}:=
\int_{\widetilde{Mp}_{0,2|n}(N,d)}ev_{0}^{*}(h^{a})\cdot ev_{\infty}^{*}(h^{b})\cdot\prod_{j=1}^{n}ev_{j}^{*}(h^{m_j}),
\label{defw}
\end{equation}
where $\cdot$ is the product of the cohomology ring $H^{*}(\widetilde{Mp}_{0,2|n}(N,d))$. We introduce a 
${\bf C}^{\times}$ action on $\widetilde{Mp}_{0,2|n}(N,d)$ to apply localization technique. First, we 
define it on the bulk stratum $Mp_{0,2|n}(N,d)$.
\begin{equation}
(e^t)\cdot [\bigl(({\bf a}_0,{\bf a}_1,\cdots,{\bf a}_d),(z_1,\cdots,z_n)\bigr)]:=
[\bigl((e^{\lambda_0 t}{\bf a}_0,e^{\lambda_1 t}{\bf a}_1,\cdots,e^{\lambda_d t}{\bf a}_d),(z_1,\cdots,z_n)\bigr)],\;(t,\lambda_i\in {\bf C}),
\label{torus1}
\end{equation}
where $\lambda_i\;(i=0,1,\cdots,d)$ are the characters of ${\bf C}^{\times}$ action.
The ${\bf C}^{\times}$ action acts only on parameters of quasi maps. The part of $\widetilde{Mp}_{0,2|n}(N,d)$
that describe parameters of quasi maps is the same as $\widetilde{Mp}_{0,2}(N,d)$ whose boundary structure is 
given by ordered partition $(d_1,\cdots,d_l)\in OP_d$. 
In \cite{mmg}, we gave toric construction of $\widetilde{Mp}_{0,2}(N,d)$ by introducing boundary divisor 
coordinates $u_j\; (j=1,2,\cdots, d-1)$. A point in $\widetilde{Mp}_{0,2}(N,d)$ was described by 
$$[({\bf a}_0,{\bf a}_1,\cdots,{\bf a}_d,u_1,\cdots,u_{d-1})]$$
where $[*]$ represents equivalence class under the $({\bf C}^{\times})^{d+1}$ action used in the toric construction.  
See \cite{mmg} for details. In \cite{mmg}, we used fundamentally the same ${\bf C}^{\times}$ action as (\ref{torus1})
, that acts on $\widetilde{Mp}_{0,2}(N,d)$. It was defined by,
\begin{equation}
(e^t)\cdot[({\bf a}_0,{\bf a}_1,\cdots,{\bf a}_d,u_1,\cdots,u_{d-1})]:=
[(e^{\lambda_0 t}{\bf a}_0,e^{\lambda_1 t}{\bf a}_1,\cdots,e^{\lambda_d t}{\bf a}_d,u_1,\cdots,u_{d-1})].
\label{torus0}
\end{equation}
Note that (\ref{torus0}) is defined on the whole strata of $\widetilde{Mp}_{0,2}(N,d)$. 
Since the part of $\widetilde{Mp}_{0,2|n}(N,d)$ describing parameters of quasi map is the same as $\widetilde{Mp}_{0,2}(N,d)$,
we can extend the ${\bf C}^{\times}$ action given by (\ref{torus1}) to whole $\widetilde{Mp}_{0,2|n}(N,d)$. 
 
Next, we determine fixed point sets of $\widetilde{Mp}_{0,2|n}(N,d)$ under the above ${\bf C}^{\times}$ action. 
We consider first the stratum $Mp_{0,2|n}(N,d)$. Since $[*]$ in (\ref{torus1}) represents equivalence class under 
$({\bf C}^{\times})^{2}$ action, we can trivialize this action on ${\bf a}_{0}$ and ${\bf a}_d$ by 
regarding them as $[{\bf a}_{0}],[{\bf a}_d]\in CP^{N-1}$. After this trivialization, the ${\bf C}^{\times}$ action
in (\ref{torus1}) is rewritten as follows:
\begin{eqnarray}
&&(e^t)\cdot \bigl(([{\bf a}_0],{\bf a}_1,\cdots,[{\bf a}_d]]),(z_1,\cdots,z_n)\bigr)=\no\\
&&\bigl(([{\bf a}_0],e^{((\lambda_1-\lambda_0)-\frac{\lambda_d-\lambda_0}{d})t}{\bf a}_1,\cdots,e^{((\lambda_{d-1}-\lambda_0)-\frac{(d-1)(\lambda_d-\lambda_0)}{d})t}{\bf a}_{d-1},[{\bf a}_d]]),(e^{(\frac{\lambda_d-\lambda_0}{d})t}z_1,\cdots,e^{(\frac{\lambda_d-\lambda_0}{d})t}z_n)\bigr)\no\\
\label{trtr}
\end{eqnarray}
Therefore, fixed points appear only if $n=0$ and they are given by,
\begin{equation}
\bigl(([{\bf a}_0],{\bf 0},\cdots,{\bf 0},[{\bf a}_d]])\bigr).
\label{fbulk}
\end{equation}
Even after the trivialization, we still have remaining ${\bf Z}_{d}$ action:
\begin{eqnarray}
&&\zeta\cdot \bigl(([{\bf a}_0],{\bf a}_1,\cdots,{\bf a}_{d-1},[{\bf a}_d]),(z_1,\cdots,z_{n})\bigr)=\no\\
&&\bigl(([{\bf a}_0],\zeta{\bf a}_1,\cdots,\zeta^{d-1}{\bf a}_{d-1},[{\bf a}_d]),(\zeta^{-1}z_1,\cdots,\zeta^{-1}z_{n})\bigr),\;\;(\zeta=\exp(\frac{2\pi\sqrt{-1}}{d})).
\end{eqnarray}
Hence, the fixed points in (\ref{fbulk}) are ${\bf Z}_d$ orbifold singularities in $\widetilde{Mp}_{0,2|0}(N,d)$. 
(\ref{fbulk}) also tells us that the fixed point set of $\widetilde{Mp}_{0,2|0}(N,d)$ is given by $CP^{N-1}\times CP^{N-1}$.
Now, we can determine the fixed point set in the stratum labeled by $(d_1,\cdots,d_l)$ and  $\bigl(\coprod_{i=0}^{l}A_i\bigr)\coprod\bigl(\coprod_{i=1}^{l}B_i\bigr)=\{1,2,\cdots,n\}$.
We have non-empty fixed point set only if $B_{i}=\emptyset\;(i=1,2,\cdots,l)$. 
If this condition is satisfied, the fixed point set of the stratum is given by $\prod_{i=0}^{l}(CP^{N-1})_{i}\times\prod_{i=0}^{l}
\overline{M}_{0,2||A_i|}$ (if $A_i=\emptyset$, we don't include $\overline{M}_{0,2||A_i|}$ in the product).  Here the 
first factor represents a chain of quasi maps:
\begin{equation}
\mathop{\cup}_{i=1}^{l}[{\bf a}_{\sum_{j=1}^{i-1}d_j}(s_i)^{d_i}+{\bf a}_{\sum_{j=1}^{i}d_j}(t_i)^{d_i}],
\end{equation}
and the second factor describes degrees of freedom of marked points that are mapped to $[{\bf a}_{\sum_{j=1}^{i}d_j}]\;(i=0,1,\cdots,l)$.
This fixed point set is also a set of orbifold singularities on which $\prod_{j=1}^{l}{\bf Z}_{d_j}$ acts. 

We have determined the fixed point sets of $\widetilde{Mp}_{0,2|n}(N,d)$ under the ${\bf C}^{\times}$ action of (\ref{torus1}).
To proceed the localization technique, we analyze contributions from normal bundle of the 
fixed point set labeled by $(d_1,\cdots,d_l)$ and $\coprod_{i=0}^{l}A_i=\{1,2,\cdots,n\}$, to localized integrand.  
For brevity, we introduce another notation of an ordered partition $(d_1,\cdots,d_l)$:
\begin{equation}
0=f_0<f_1<f_2<\cdots<f_{l-1}<f_l=d,\;\;f_j-f_{j-1}=d_j,\;\;(j=1,2\cdots,l).
\label{aop}
\end{equation}
In the following, we denote by $h_{f_i}$ the hyperplane class of $(CP^{N-1})_i$ in 
$\prod_{i=0}^{l}(CP^{N-1})_{i}\times\prod_{i=0}^{l}\overline{M}_{0,2||A_i|}$.

We first compute contribution from $Mp_{0,2|0}(N,f_j-f_{j-1})$ in (\ref{mpdn}). By fixing the ambiguity coming 
from $({\bf C}^{\times})^{2}$ action, we can represent $Mp_{0,2|0}(N,f_j-f_{j-1})$ in the following form:
\begin{equation}
Mp_{0,2|0}(N,f_j-f_{j-1})=\{([{\bf a}_{f_{i-1}}],{\bf y}_{f_{i-1}+1},{\bf y}_{f_{i-1}+2},\cdots,{\bf y}_{f_{i}-1},[{\bf a}_{f_{i}}])\;
|\;[{\bf a}_{f_{i-1}}],[{\bf a}_{f_{i}}]\in CP^{N-1},\;{\bf y}_j\in {\bf C}^{N}\}/{\bf Z}_{f_{j}-f_{j-1}}.
\label{noamp}
\end{equation}
Therefore, normal bundle is given by $\displaystyle{\mathop{\oplus}_{i=1}^{f_j-f_{j-1}-1} \mathop{\oplus}_{l=1}^{N}\frac{\d}{\d y_{f_{j-1}+i}^{l}}}$.
As we have discussed in \cite{mmg}, $\displaystyle{\frac{\d}{\d y_{f_{j-1}+i}^{l}}}$ is isomorphic to 
${\cal O}_{(CP^{N-1})_{j-1}}(\frac{f_{j}-f_{j-1}-i}{f_{j}-f_{j-1}})\otimes {\cal O}_{(CP^{N-1})_{j}}(\frac{i}{f_{j}-f_{j-1}}) $ as an 
orbi-bundle on $Mp_{0,2|0}(N,f_j-f_{j-1})$. Hence its first Chern class is given by,
\begin{equation}
\bigl( \frac{f_{j}-f_{j-1}-i}{f_{j}-f_{j-1}}\bigr)h_{f_{j-1}}+\bigl( \frac{i}{f_{j}-f_{j-1}}\bigr)h_{f_{j}}.
\end{equation}
According to \cite{mmg}, its character of the ${\bf C}^{\times}$ action is, 
\begin{equation}
\bigl( \frac{f_{j}-f_{j-1}-i}{f_{j}-f_{j-1}}\bigr)\lambda_{f_{j-1}}+\bigl( \frac{i}{f_{j}-f_{j-1}}\bigr)\lambda_{f_{j}}-\lambda_{f_{j-1}+i}.
\end{equation} 
These results lead us to the following contribution to the localized integrand:
\begin{equation}
\frac{1}{\displaystyle{\prod_{i=1}^{f_{j}-f_{j-1}-1}\biggl
(\bigl( \frac{f_{j}-f_{j-1}-i}{f_{j}-f_{j-1}}\bigr)(h_{j-1}+\lambda_{f_{j-1}})+\bigl( \frac{i}{f_{j}-f_{j-1}}\bigr)(h_{j}+\lambda_{f_{j}})-\lambda_{f_{j-1}+i}
\biggr)^{N}}}.
\label{cont1}
\end{equation}
Next, we compute contribution from $\overline{M}_{0,2||A_i|},\;(i=1,2,\cdots,l-1)$.
If $A_i=\emptyset$, the contribution comes from smoothing nodal singularity $[{\bf a}_{f_i}]$. 
This factor in the normal bundle is identified with $\frac{d}{d(\frac{t_i}{s_{i}})}\otimes
\frac{d}{d(\frac{s_{i+1}}{t_{i+1}})}$ and its equivariant first Chern class is 
$\frac{h_{f_i}+\lambda_{f_i}-h_{f_{i-1}}-\lambda_{f_{i-1}}}{f_{i}-f_{i-1}}+
\frac{h_{f_i}+\lambda_{f_i}-h_{f_{i+1}}-\lambda_{f_{i+1}}}{f_{i+1}-f_{i}}$. Therefore, 
the contribution to the localized integrand is given as follows:
\begin{equation}
\frac{1}{\displaystyle{\frac{h_{f_i}+\lambda_{f_i}-h_{f_{i-1}}-\lambda_{f_{i-1}}}{f_{i}-f_{i-1}}+
\frac{h_{f_i}+\lambda_{f_i}-h_{f_{i+1}}-\lambda_{f_{i+1}}}{f_{i+1}-f_{i}}}}.
\label{cont2}
\end{equation}
If $A_i\neq \emptyset$, we have a nontrivial stable curve that corresponds to a point in $\overline{M}_{0,2||A_i|}$.
The factors in the normal bundle coming from these components are degrees of freedom of smoothing nodal 
singularities that connect the stable curve with  
$[{\bf a}_{f_{i-1}}(s_i)^{f_i-f_{i-1}}+{\bf a}_{f_{i}}(t_i)^{f_i-f_{i-1}}]$ and $[{\bf a}_{f_{i}}(s_{i+1})^{f_{i+1}-f_{i}}+{\bf a}_{f_{i+1}}(t_{i+1})^{f_{i+1}-f_{i}}]$.  
Let $C_i$ be the stable curve mentioned above. Then these two factors are identified with $\frac{d}{d(\frac{t_i}{s_{i}})}\otimes
T^{\prime}_{0}C_i$ and $T^{\prime}_{\infty}C_i\otimes
\frac{d}{d(\frac{s_{i+1}}{t_{i+1}})}$. Hence we obtain the following contribution:
\begin{equation}
\frac{1}{\displaystyle{\biggl(\frac{h_{f_i}+\lambda_{f_i}-h_{f_{i-1}}-\lambda_{f_{i-1}}}{f_{i}-f_{i-1}}+
c_{1}(T^{\prime}_{0}C_i)\biggr)
\biggl(\frac{h_{f_i}+\lambda_{f_i}-h_{f_{i+1}}-\lambda_{f_{i+1}}}{f_{i+1}-f_{i}}+c_{1}(T^{\prime}_{\infty}C_i)\biggr)}}.
\label{cont3p}
\end{equation}
At this stage, we temporarily set $\Lambda_{i0}:=\frac{h_{f_i}+\lambda_{f_i}-h_{f_{i-1}}-\lambda_{f_{i-1}}}{f_{i}-f_{i-1}},\;\;
\Lambda_{i\infty}=\frac{h_{f_i}+\lambda_{f_i}-h_{f_{i+1}}-\lambda_{f_{i+1}}}{f_{i+1}-f_{i}}$ and rewrite (\ref{cont3p}) as follows:
\begin{equation}
\frac{1}{(\Lambda_{i0}-
c_{1}(T^{\prime *}_{0}C_i))(
\Lambda_{i\infty}-c_{1}(T^{\prime *}_{\infty}C_i))}.
\label{cont3r}
\end{equation}
To proceed the localization technique, we integrate out the above equivariant class on $\overline{M}_{0,2||A_i|}$.
\begin{equation}
\int_{\overline{M}_{0,2||A_i|}}\frac{1}{(\Lambda_{i0}-
c_{1}(T^{\prime *}_{0}C_i))(
\Lambda_{i\infty}-c_{1}(T^{\prime *}_{\infty}C_i))}=\frac{1}{\Lambda_{i0}\Lambda_{i\infty}}
\cdot\sum_{n,m=0}^{\infty}\frac{1}{(\Lambda_{i0})^n(\Lambda_{i\infty})^m}\int_{\overline{M}_{0,2||A_i|}}
(c_{1}(T^{\prime *}_{0}C_i))^n(c_{1}(T^{\prime *}_{\infty}C_i))^m.
\label{medcomp}
\end{equation}  
The intersection number that appear at the right end of (\ref{medcomp}) has been already computed in 
\cite{alex}, \cite{pand} and it is given as follows:
\begin{equation}
\int_{\overline{M}_{0,2||A_i|}}
(c_{1}(T^{\prime *}_{0}C_i))^n(c_{1}(T^{\prime *}_{\infty}C_i))^m=\delta_{|A_i|-1,n+m}\cdot{n+m\choose n}.
\label{int}
\end{equation}
Combining (\ref{medcomp}) with (\ref{int}), we obtain the following contribution to the localized integrand:
\begin{eqnarray} 
&&\frac{1}{\Lambda_{i0}\Lambda_{i\infty}}(\frac{1}{\Lambda_{i0}}+\frac{1}{\Lambda_{i\infty}})^{|A_i|-1}=
\frac{1}{\displaystyle{\frac{h_{f_i}+\lambda_{f_i}-h_{f_{i-1}}-\lambda_{f_{i-1}}}{f_{i}-f_{i-1}}+
\frac{h_{f_i}+\lambda_{f_i}-h_{f_{i+1}}-\lambda_{f_{i+1}}}{f_{i+1}-f_{i}}}}\times\no\\
&&\biggl(\displaystyle{\frac{f_{i}-f_{i-1}}{h_{f_i}+\lambda_{f_i}-h_{f_{i-1}}-\lambda_{f_{i-1}}}+
\frac{f_{i+1}-f_{i}}{h_{f_i}+\lambda_{f_i}-h_{f_{i+1}}-\lambda_{f_{i+1}}}}\biggr)^{|A_i|}.
\label{cont3}
\end{eqnarray}
If we compare (\ref{cont3}) with (\ref{cont2}), we can easily see that (\ref{cont3}) is applicable to the case of $A_i=\emptyset$.
We also have to determine contributions from $\overline{M}_{0,2||A_{0}|}$ and $\overline{M}_{0,2||A_{l}|}$.
It is sufficient to consider the case of $\overline{M}_{0,2||A_{0}|}$. If $A_{0}=\emptyset$, we have no 
contribution from the normal bundle to the localized integrand. If $A_{0}\neq\emptyset$, contribution 
comes from smoothing the nodal singularity that connects the stable curve described by $\overline{M}_{0,2||A_{0}|}$
with $[{\bf a}_{f_{0}}(s_1)^{f_1-f_{0}}+{\bf a}_{f_{1}}(t_i)^{f_1-f_{0}}]$. As discussed in the previous case, we have to 
compute the following integral:
\begin{equation}
\int_{\overline{M}_{0,2||A_{0}|}}\frac{1}{\displaystyle{\biggl(\frac{h_{f_0}+\lambda_{f_0}-h_{f_{1}}-\lambda_{f_{1}}}{f_{1}-f_{0}}+c_{1}(T^{\prime}_{\infty}C_1)\biggr)}}.
\label{cont4p}
\end{equation}
Application of (\ref{int}) leads us to the following contribution to the localized integrand:
\begin{equation}
\displaystyle{\biggl(\frac{f_{1}-f_{0}}{h_{f_0}+\lambda_{f_0}-h_{f_{1}}-\lambda_{f_{1}}}\biggr)^{|A_{0}|}}.
\label{cont4i}
\end{equation}
In the same way as above, we obtain the contribution from $\overline{M}_{0,2||A_{l}|}$:
\begin{equation}
\displaystyle{\biggl(\frac{f_{l}-f_{l-1}}{h_{f_l}+\lambda_{f_l}-h_{f_{l-1}}-\lambda_{f_{l-1}}}\biggr)^{|A_{l}|}}.
\label{cont4f}
\end{equation}
Finally, the contributions to localized integrand coming from $ev_{0}(h^a)$, $ev_{\infty}(h^b)$ and $ev_{i}(h^{m_i}),\;(i=1,2,\cdots,n)$
are given by the following correspondence:
\begin{eqnarray}
&&ev_{0}(h^a)\mapsto (h_{f_{0}}+\lambda_{f_{0}})^a,\no\\
&&ev_{\infty}(h^b)\mapsto (h_{f_{l}}+\lambda_{f_{l}})^b, \no\\
&&ev_{i}(h^{m_i})\mapsto (h_{f_{j}}+\lambda_{f_{j}})^{m_i},
\label{loev}
\end{eqnarray}    
where we assume that the marked point $z_i$ is mapped to $[{\bf a}_{f_{j}}]\in CP^{N-1}$. 

What remains to complete localization computation is to combine the factors given in (\ref{cont1}) ,(\ref{cont3}), (\ref{cont4i}), (\ref{cont4f}) and (\ref{loev}) 
into the localized integrand labeled by $(d_1,\cdots,d_l)$ and $\coprod_{j=0}^{l}A_j=\{1,2,\cdots,n\}$, to integrate it out on $\prod_{j=1}^{l}(CP^{N-1})_j$
 and to sum up the results of integration by labels. We now discuss the label $\coprod_{j=0}^{l}A_j=\{1,2,\cdots,n\}$. Since $A_j$ can be an empty set,
 this label is equivalent to determining the point $[{\bf a}_{f_j}]$ to which the marked point $z_i$ is mapped, for each $i\in\{1,2,\cdots,n\}$.
 Therefore, we can replace the label $\coprod_{j=0}^{l}A_j=\{1,2,\cdots,n\}$ by a sequence $(j_1,j_2,\cdots,j_n),\;(j_i\in\{0,1,\cdots,l\})$.
 From this point of view, $|A_j|$ is just the number of $i$'s that satisfy the condition $j_i=j$. Hence the localized integrand labeled by 
$(d_1,\cdots,d_l)$ and $(j_1,j_2,\cdots,j_n)$ is written down as follows:
\begin{eqnarray}
&&\prod_{j=1}^{l}\frac{1}{\displaystyle{\prod_{i=1}^{d_j-1}\biggl
(\bigl(\frac{d_{j}-i}{d_{j}}\bigr)(h_{f_{j-1}}+\lambda_{f_{j-1}})+\bigl( \frac{i}{d_{j}}\bigr)(h_{f_{j}}+\lambda_{f_{j}})-\lambda_{f_{j-1}+i}
\biggr)^{N}}}\times\no\\
&&\prod_{j=1}^{l-1}\frac{1}{\displaystyle{\frac{h_{f_i}+\lambda_{f_i}-h_{f_{i-1}}-\lambda_{f_{i-1}}}{d_{i}}+
\frac{h_{f_i}+\lambda_{f_i}-h_{f_{i+1}}-\lambda_{f_{i+1}}}{d_{i+1}}}}\times\no\\
&&\prod_{i=1}^{n}\bigl(\displaystyle{\frac{d_{j_i}}{h_{f_{j_i}}+\lambda_{f_{j_i}}-h_{f_{{j_i}-1}}-\lambda_{f_{j_i-1}}}+
\frac{d_{j_i+1}}{h_{f_{j_i}}+\lambda_{f_{j_i}}-h_{f_{j_i+1}}-\lambda_{f_{j_i+1}}}}\bigr)(h_{f_{j_i}}+\lambda_{f_{j_i}})^{m_i}\times\no\\
&&(h_{0}+\lambda_0)^a(h_d+\lambda_d)^b.
\label{loint}
\end{eqnarray} 
In (\ref{loint}), we used both $(d_1,\cdots,d_l)$ and $0=f_0<\cdots<f_l=d$ to denote an ordered partition. We also formally set $d_0=d_{l+1}=0$ for brevity. 
Before integrate out (\ref{loint}), we sum up the integrand by the label $(j_1,j_2,\cdots,j_n)$.
If we pay attention to the third line of (\ref{loint}), the factor that comes from $i$ sums up to the following factor by varying
$j_i$ from $0$ to $l$:
\begin{equation}
\sum_{j=0}^{l}\bigl(\displaystyle{\frac{d_{j}}{h_{f_{j}}+\lambda_{f_{j}}-h_{f_{{j}-1}}-\lambda_{f_{j-1}}}+
\frac{d_{j+1}}{h_{f_{j}}+\lambda_{f_{j}}-h_{f_{j+1}}-\lambda_{f_{j+1}}}}\bigr)(h_{f_{j}}+\lambda_{f_{j}})^{m_i}
\label{isum}  
\end{equation} 
We introduce here a rational function:
\begin{equation}
w_{a}^{d}(z,w):=d\cdot\frac{z^a-w^a}{z-w}. 
\label{wad}
\end{equation} 
Then we can rewrite (\ref{isum}) into the form,
\begin{equation}
\sum_{j=1}^{l}w_{m_i}^{d_j}(h_{f_{j-1}}+\lambda_{f_{j-1}},h_{f_j}+\lambda_{f_j}).
\end{equation} 
These consideration leads us to the following formula for the localized integrand summed up by the label $(j_1,j_2,\cdots,j_n)$.
\begin{eqnarray}
&&(h_{0}+\lambda_0)^a(h_d+\lambda_d)^b\prod_{j=1}^{l}\frac{1}{\displaystyle{\prod_{i=1}^{d_j-1}\biggl
(\bigl(\frac{d_{j}-i}{d_{j}}\bigr)(h_{f_{j-1}}+\lambda_{f_{j-1}})+\bigl( \frac{i}{d_{j}}\bigr)(h_{f_{j}}+\lambda_{f_{j}})-\lambda_{f_{j-1}+i}
\biggr)^{N}}}\times\no\\
&&\prod_{j=1}^{l-1}\frac{1}{\displaystyle{\frac{h_{f_j}+\lambda_{f_j}-h_{f_{j-1}}-\lambda_{f_{j-1}}}{d_{j}}+
\frac{h_{f_j}+\lambda_{f_j}-h_{f_{j+1}}-\lambda_{f_{j+1}}}{d_{j+1}}}}\cdot\prod_{i=1}^{n}\bigl(\sum_{j=1}^{l}w_{m_i}^{d_j}(h_{f_{j-1}}+\lambda_{f_{j-1}},h_{f_j}+\lambda_{f_j})\bigr).\no\\
\label{loints}
\end{eqnarray} 
Final step of localization computation is to integrate the equivarinat form (\ref{loints}) and to sum up the results by the label
$(d_1,\cdots,d_l)$. We have one subtle remark here. As we have pointed out in determining fixed point sets, the fixed point set 
labeled by $(d_1,d_2,\cdots,d_l)$ is the set of orbifold singularities on which $\prod_{j=1}^{l}{\bf Z}_{d_j}$ acts. Therefore, 
we have to divide the integral that comes from $(d_1,d_2,\cdots,d_l)$, by the factor $\prod_{j=1}^{l}d_j$ and we obtain, 
\begin{eqnarray}
&&w({\cal O}_{h^a}{\cal O}_{h^b}|\prod_{i=1}^{n}{\cal O}_{h^{m_i}})_{0,d}=\sum_{(d_1,\cdots,d_l)\in OP_d}\biggl(\prod_{j=1}^{l}\frac{1}{d_j}\biggr)
\int_{(CP^{N-1})_0}\int_{(CP^{N-1})_1}\cdots\int_{(CP^{N-1})_l}\times\no\\
&&(h_{0}+\lambda_0)^a(h_d+\lambda_d)^b\prod_{j=1}^{l}\frac{1}{\displaystyle{\prod_{i=1}^{d_j-1}\biggl
(\bigl(\frac{d_{j}-i}{d_{j}}\bigr)(h_{f_{j-1}}+\lambda_{f_{j-1}})+\bigl( \frac{i}{d_{j}}\bigr)(h_{f_{j}}+\lambda_{f_{j}})-\lambda_{f_{j-1}+i}
\biggr)^{N}}}\times\no\\
&&\prod_{j=1}^{l-1}\frac{1}{\displaystyle{\frac{h_{f_j}+\lambda_{f_j}-h_{f_{j-1}}-\lambda_{f_{j-1}}}{d_{j}}+
\frac{h_{f_j}+\lambda_{f_j}-h_{f_{j+1}}-\lambda_{f_{j+1}}}{d_{j+1}}}}\cdot\prod_{i=1}^{n}\bigl(\sum_{j=1}^{l}w_{m_i}^{d_j}(h_{f_{j-1}}+\lambda_{f_{j-1}},h_{f_j}+\lambda_{f_j})\bigr),\no\\
\label{wints2}
\end{eqnarray} 
 where the hyperplane class of $(CP^{N-1})_{j},\;(j=0,1,\cdots,l)$ is given by $h_{f_j}$.
 
 But this is not the end of the story. If we pay attention to the equality:$\int_{CP^{N-1}}h^j=\delta_{j,N-1}$, we can apply replacements,
\begin{eqnarray}
h_{f_j}&\rightarrow& z_{f_j},\no\\
\int_{(CP^{N-1})_j}&\rightarrow& \frac{1}{2\pi\sqrt{-1}}\oint_{C_{(0)}}\frac{dz_{f_j}}{(z_{f_j})^N},
\end{eqnarray}
where $z_{f_j}$ is a complex variable and $\frac{1}{2\pi\sqrt{-1}}\oint_{C_{(0)}}dz$ is the operation of taking a residue at $z=0$. 
We then apply shift of variables $z_{f_j}\rightarrow z_{f_j}-\lambda_{f_j},\;(j=0,1,\cdots,l)$ and obtain the following equality:
\begin{eqnarray}
&&w({\cal O}_{h^a}{\cal O}_{h^b}|\prod_{i=1}^{n}{\cal O}_{h^{m_i}})_{0,d}=\no\\
&&\sum_{0=f_0<f_1<\cdots<f_l=d}\biggl(\prod_{j=1}^{l}\frac{1}{f_j-f_{j-1}}\biggr)\frac{1}{(2\pi\sqrt{-1})^{l+1}}
\oint_{C_{(\lambda_{f_0})}}\frac{dz_{f_0}}{(z_{f_0}-\lambda_{f_0})^N} 
\cdots 
\oint_{C_{(\lambda_{f_l})}}\frac{dz_{f_l}}{(z_{f_l}-\lambda_{f_l})^N}
\times\no\\
&&(z_{0})^a(z_{d})^b\prod_{j=1}^{l}\frac{1}{\displaystyle{\prod_{i=1}^{d_j-1}\biggl
(\bigl(\frac{d_{j}-i}{d_{j}}\bigr)z_{f_{j-1}}+\bigl(\frac{i}{d_{j}}\bigr)z_{f_{j}}-\lambda_{f_{j-1}+i}
\biggr)^{N}}}\prod_{j=1}^{l-1}\frac{1}{\displaystyle{\frac{z_{f_j}-z_{f_{j-1}}}{d_{j}}+
\frac{z_{f_j}-z_{f_{j+1}}}{d_{j+1}}}}\times\no\\
&&\prod_{i=1}^{n}\bigl(\sum_{j=1}^{l}w_{m_i}^{d_j}(z_{f_{j-1}},z_{f_j})\bigr),
\label{wres1}
\end{eqnarray} 
where $\frac{1}{2\pi\sqrt{-1}}\oint_{C_{(\lambda)}}dz$ means the operation of taking a residue at $z=\lambda$.
Let us consider here the following residue integral:
\begin{eqnarray}
&&\frac{1}{(2\pi\sqrt{-1})^{d+1}}\oint_{E^0_{(\lambda_{0})}}\frac{dz_0}{(z_0-\lambda_0)^N}\oint_{E^1_{(\lambda_1)}}\frac{dz_1}{(z_1-\lambda_1)^N}\cdots \oint_{E^d_{(\lambda_d)}}\frac{dz_d}{(z_d-\lambda_d)^N}
\times\no\\
&&(z_{0})^a\cdot\biggl(\prod_{j=1}^{d-1}\frac{1}{(2z_j-z_{j-1}-z_{j+1})}\biggr)\cdot(z_{d})^{b}\cdot
\prod_{i=1}^{n}\biggl(\sum_{j=1}^{d}w_{m_i}^1(z_{j-1},z_j)\biggr),
\label{key}
\end{eqnarray}
where $\frac{1}{2\pi\sqrt{-1}}\oint_{E^j_{(\lambda_j)}}dz_j$ means the operation of taking residues at $z_j=\lambda_j$ and $z_j=\frac{z_{j-1}+z_{j+1}}{2}$ for 
$j=1,2,\cdots,l-1$ (resp. $z_j=\lambda_j$ for $j=0,d$). As we demonstrated in \cite{mmg}, we can observe by elementary computation that 
the summand labeled by $0=f_0<f_1<\cdots<f_l=d$ in (\ref{wres1}) is obtained by taking residues of (\ref{key}) at $z_{j}=\frac{z_{j-1}+z_{j+1}}{2}\;(j\in \{0,1,\cdots,d\}\setminus
\{f_0,f_1,\cdots,f_l\})$ and at $z_{j}=\lambda_j\;(j\in\{f_0,f_1,\cdots,f_l\})$ respectively. Hence (\ref{key}) equals 
$w({\cal O}_{h^a}{\cal O}_{h^b}|\prod_{i=1}^{n}{\cal O}_{h^{m_i}})_{0,d}$. Finally, we take the non-equivariant limit $\lambda_j\rightarrow 0\;(j=0,1,\cdots,d)$
of (\ref{key}) 
and obtain the following theorem:  
\begin{theorem}
\begin{eqnarray}
w({\cal O}_{h^a}{\cal O}_{h^{b}}|\prod_{i=1}^{n}({\cal O}_{h^{m_i}}))_{0,d}&=&
\frac{1}{(2\pi\sqrt{-1})^{d+1}}\oint_{E^0_{(0)}}\frac{dz_0}{(z_0)^N}\oint_{E^1_{(0)}}\frac{dz_1}{(z_1)^N}\cdots \oint_{E^d_{(0)}}\frac{dz_d}{(z_d)^N}
\times\no\\
&&(z_{0})^a\cdot\biggl(\prod_{j=1}^{d-1}\frac{1}{(2z_j-z_{j-1}-z_{j+1})}\biggr)\cdot(z_{d})^{b}\cdot
\prod_{i=1}^{n}\biggl(\sum_{j=1}^{d}w_{m_i}^1(z_{j-1},z_j)\biggr), \;\;(d>0).\no\\
\label{th1}
\end{eqnarray}
\end{theorem}     
 
\subsection{Numerical Computation in the $CP^2$ case}
In this section, we discuss how to compute the genus $0$ Gromov-Witten invariants of $CP^2$ by using the multi-point virtual structure 
constants, i.e., the intersection 
numbers computed in the previous section. $H^{*,*}(CP^2)$ is spanned by $h^j$ $(j=0,1,2)$. Therefore it is convenient
to write the intersection number as $w({\cal O}_{h^a}{\cal O}_{h^{b}}|\prod_{j=0}^{2}({\cal O}_{h^{j}})^{m_j})_{0,d}$.
The formula (\ref{th1}) restricted to $CP^2$ is given as follows.
\begin{eqnarray}
w({\cal O}_{h^a}{\cal O}_{h^{b}}|\prod_{j=0}^{2}({\cal O}_{h^{j}})^{m_j})_{0,d}&=&
\frac{1}{(2\pi\sqrt{-1})^{d+1}}\oint_{E^0_{(0)}}\frac{dz_0}{(z_0)^3}\oint_{E^1_{(0)}}\frac{dz_1}{(z_1)^3}\cdots \oint_{E^d_{(0)}}\frac{dz_d}{(z_d)^3}
\times\no\\
&&(z_{0})^a\cdot\biggl(\prod_{j=1}^{d-1}\frac{1}{(2z_j-z_{j-1}-z_{j+1})}\biggr)\cdot(z_{d})^{b}\cdot
\prod_{j=0}^{2}\biggl(\sum_{i=1}^{d}w_{j}^1(z_{i-1},z_i)\biggr)^{m_j}, \;\;(d>0).\no\\
\label{p2c}
\end{eqnarray}
We show below important characteristics of the intersection number that follows from (\ref{p2c}).
\begin{prop} The multi-point virtual structure constants satisfy the following equalities:\\ 
(i)
\begin{equation}
 w({\cal O}_{h^a}{\cal O}_{h^{b}}|{\cal O}_{1}\prod_{j=0}^{2}({\cal O}_{h^{j}})^{m_j})_{0,d}=0.
\label{punc}
\end{equation}
(ii)
\begin{equation}
 w({\cal O}_{h^a}{\cal O}_{h^{b}}|{\cal O}_{h}\prod_{j=0}^{2}({\cal O}_{h^{j}})^{m_j})_{0,d}=
d\cdot w({\cal O}_{h^a}{\cal O}_{h^{b}}|\prod_{j=0}^{2}({\cal O}_{h^{j}})^{m_j})_{0,d}.
\label{divisor}
\end{equation}
\label{divpun}
\end{prop}
{\it proof)} (i) follows from the equality :$\displaystyle{w_{0}^{1}(z,w)=\frac{1-1}{z-w}=0}$. As for (ii), it is enough to note,
$$ \sum_{i=1}^{d}w_{1}^{1}(z_{i-1},z_{i})=\sum_{j=1}\frac{z_{i-1}-z_{i}}{z_{i-1}-z_{i}}=d.\;\;\;\;\;\;\; \Box$$
This proposition says that for operator insertions at the right hand side of "$|$" , both the divisor axiom and the puncture axiom for Gromov-Witten 
invariants hold. Motivated by this fact, we introduce the generating function of the multi-point virtual structure constants 
with respect to operator insertions at the right hand side of "$|$".
\begin{defi}
\begin{eqnarray}
w({\cal O}_{h^a}{\cal O}_{h^{b}}|(x^0,x^1,x^2))_{0}&:=&x^c\cdot\int_{CP^2}h^{a+b+c}+
\sum_{d>0,\{m_j\}}w({\cal O}_{h^a}{\cal O}_{h^{b}}|\prod_{j=0}^{2}({\cal O}_{h^{j}})^{m_j})_{0,d}\cdot\prod_{j=0}^{2}\frac{(x^{j})^{m_j}}{m_j!},
\label{genc}
\end{eqnarray}
where $x^{j}$ $(j=0,1,2)$ is the variable associated with insertion of ${\cal O}_{h^j}$. 
\end{defi}
In (\ref{genc}), we only consider three operator insertions for degree $0$ virtual structure constants and identify them with classical 
intersection numbers, just as we do in the case of Gromov-Witten invariants. The assertion of Proposition 1 simplifies the generating function 
in the same way as the Gromov-Witten case. 
\begin{prop}
\begin{equation}
w({\cal O}_{h^a}{\cal O}_{h^{b}}|(x^0,x^1,x^2))_{0}=x^c\cdot\int_{CP^2}h^{a+b+c}+\sum_{d>0,m_2}
w({\cal O}_{h^a}{\cal O}_{h^{b}}|({\cal O}_{h^{2}})^{m})_{0,d}\cdot e^{dx^1}\cdot\frac{(x^{2})^{m}}{m!}.
\end{equation}
\end{prop}
Therefore, we only have to compute $w({\cal O}_{h^a}{\cal O}_{h^{b}}|({\cal O}_{h^{2}})^{m})_{0,d}$. From degree counting of the residue 
integral formula, we can easily see that the non-zero virtual structure constants appear only when $a+b+m-2=3d-1$ for $d\geq 1$.
We also introduce here the corresponding generating functions of genus $0$ Gromov-Witten invariant of $CP^2$. 
\begin{defi}
Let $\langle\prod_{j=0}^{2}({\cal O}_{h^j})^{m_j}\rangle_{0,d}$ be the rational Gromov-Witten invariant of degree $d$ of $CP^{2}$.
\begin{eqnarray}
\langle{\cal O}_{h^a}{\cal O}_{h^{b}}(t^0,t^1,t^2)\rangle_{0}&:=&t^c\cdot\int_{CP^2}h^{a+b+c}+
\sum_{d>0,\{m_j\}}\langle{\cal O}_{h^a}{\cal O}_{h^{b}}\prod_{j=0}^{2}({\cal O}_{h^{j}})^{m_j}\rangle_{0,d}\cdot\prod_{j=0}^{2}\frac{(t^{j})^{m_j}}{m_j!}\no\\
&=&t^c\cdot\int_{CP^2}h^{a+b+c}+\sum_{d>0,m}
\langle{\cal O}_{h^a}{\cal O}_{h^{b}}({\cal O}_{h^{2}})^{m}\rangle_{0,d}\cdot e^{dt^1}\cdot\frac{(t^{2})^{m}}{m!},
\end{eqnarray}
where $t^{j}$ $(j=0,1,2)$ is the variable associated with insertion of ${\cal O}_{h^j}$. 
\end{defi}
If we use the usual generating function of genus Gromov-Witten invariants of $CP^2$ : $F_{CP^2}(t^0,t^1,t^2)$, the  above generating
function $\langle{\cal O}_{h^a}{\cal O}_{h^{b}}(t^0,t^1,t^2)\rangle_{0}$ is nothing but $\displaystyle{\frac{\d^2 F_{CP^2}}{\d t^a \d t^b}}$. 
Therefore, we obtain integrable condition:
\begin{equation}
\frac{\d}{\d t^a}\langle{\cal O}_{h^b}{\cal O}_{h^{c}}(t^0,t^1,t^2)\rangle_{0}=\frac{\d}{\d t^b}\langle{\cal O}_{h^a}{\cal O}_{h^{c}}(t^0,t^1,t^2)\rangle_{0}.
\end{equation}
Our question in this section is whether we can compute $\langle{\cal O}_{h^a}{\cal O}_{h^{b}}(t^0,t^1,t^2)\rangle_{0}$ by using 
$w({\cal O}_{h^a}{\cal O}_{h^{b}}|(x^0,x^1,x^2))_{0}$ as the starting point.
Our key to answer this question is the observation done in \cite{mmg}.
Let us illustrate our observation by taking quintic Calabi-Yau hypersurface in $CP^4$ as an example.
In \cite{mmg}, we consider the two point virtual structure constant $w({\cal O}_{h^a}{\cal O}_{h^b})_{0,d}$ for 
the hypersurface that is given as the intersection number on $\widetilde{Mp}_{0,2}(5,d)$. We introduced 
the generating function, 
\begin{equation}
w({\cal O}_{h^a}{\cal O}_{h^b})_{0}(x):=5x\cdot \delta_{a+b,2}+\sum_{d=1}^{\infty}w({\cal O}_{h^a}{\cal O}_{h^b})_{0,d}e^x,
\end{equation}
which is non-zero only when $a+b=2$.
We proved  that $t=\frac{1}{5}w({\cal O}_{h^2}{\cal O}_{1})_{0}(x)$ coincides with the mirror map used in the mirror computation 
of Gromov-Witten invariants of the hypersurface and that $w({\cal O}_{h}{\cal O}_{h})_{0}(x(t))$ gives the generating 
function of Gromov-Witten invariants $5t+\sum_{d=1}^{\infty}\langle{\cal O}_{h}{\cal O}_{h}\rangle_{0,d}e^{dt}$.
We introduce here the classical intersection intersection matrix $\eta_{ab}=\int_{CP^4}5h^{a+b+1}=5\delta_{a+b,3}$ and 
its inverse $\eta^{ab}=\frac{1}{5}\delta_{a+b,3}$. By identifying $t$ (resp. $x$) with $t^1$ (resp. $x^1$), we conjectured
generally that 
$$ t^a=\sum_{b}\eta^{ab}w({\cal O}_{h^b}{\cal O}_{1})_{0}(x^1) $$
gives the mirror map for the mirror computation.  Since $CP^2$ is a Fano manifold with $c_1(CP^2)=3h$, $w({\cal O}_{h^b}{\cal O}_{1})_{0,d}=0$ 
for $d\geq 1$. Therefore, we have trivial mirror map if we only consider the two point virtual structure constants. 
That's why we introduce the multi-point virtual structure constants. If we pay attention to the fact that 
the classical intersection matrix of $CP^2$ is given by $\eta^{ab}=\eta_{ab}=\delta_{a+b,2}$, we are naturally 
led to propose the following conjecture. 
\begin{conj}
If we define the mirror map,
\begin{eqnarray}
t^{j}(x^0,x^1,x^2):=w({\cal O}_{h^{2-j}}{\cal O}_{1}|(x^0,x^1,x^2))_{0},
\end{eqnarray} 
we have the following equality:
\begin{equation}
\langle{\cal O}_{h^a}{\cal O}_{h^{b}}(t^0(x^0,x^1,x^2),t^1(x^0,x^1,x^2),t^2(x^0,x^1,x^2))\rangle_{0}
=w({\cal O}_{h^a}{\cal O}_{h^{b}}|(x^0,x^1,x^2))_{0}.
\end{equation}
Conversely, if we invert the mirror map,
\begin{equation}
x^{j}=x^{j}(t^0,t^1,t^2),
\end{equation}
we obtain the mirror formula to compute the rational Gromov-Witten invariants of $CP^{2}$ from the multi-point virtual structure 
constants:
\begin{equation}
\langle{\cal O}_{h^a}{\cal O}_{h^{b}}(t^0,t^1,t^2)\rangle_{0}
=w({\cal O}_{h^a}{\cal O}_{h^{b}}|(x^0(t^0,t^1,t^2),x^1(t^0,t^1,t^2),x^2(t^0,t^1,t^2)))_{0}.
\end{equation}
\end{conj}
By using the formula (\ref{p2c}), we obtain the mirror maps explicitly.
\begin{eqnarray}
{t^{2}} &=& {x^{2}} + {\displaystyle \frac {1}{4}} q(x^{2})^{4
} + {\displaystyle \frac {33}{70}} q^{2}(x^{2})^{7} + 
{\displaystyle \frac {16589}{12600}} q^{3}(x^{2})^{10} + 
{\displaystyle \frac {143698921}{32432400}} q^{4}(x^{2})^{13}
 + {\displaystyle \frac {75631936691}{4540536000}} q^{5}(x^{2
})^{16}+\cdots,\no\\
{t^{1}} &=& {x^{1}} + {\displaystyle \frac {1}{2}} (x^{2})^{3}
q + {\displaystyle \frac {7}{10}} (x^{2})^{6}q^{2} + 
{\displaystyle \frac {2593}{1512}} q^{3}(x^{2})^{9} + 
{\displaystyle \frac {2668063}{498960}} q^{4}(x^{2})^{12} + 
{\displaystyle \frac {120501923}{6306300}} q^{5}(x^{2})^{15}+\cdots,\no\\
{t^{0}} &=& {x^{0}} + {\displaystyle \frac {1}{2}} (x^{2})^{2}
q + {\displaystyle \frac {8}{15}} (x^{2})^{5}q^{2} + 
{\displaystyle \frac {983}{840}} q^{3}(x^{2})^{8} + 
{\displaystyle \frac {4283071}{1247400}} q^{4}(x^{2})^{11} + 
{\displaystyle \frac {4019248213}{340540200}} q^{5}(x^{2})^{
14}+\cdots,\no\\
&&\hspace{12cm}(q:=e^{x^1}).
\label{cp2mirr}
\end{eqnarray}
Of course, we can also compute one of the generating function,
\begin{eqnarray}
&&w({\cal O}_{h}{\cal O}_{h}|(x^0,x^1,x^2))_{0}=\no\\
&&=x^{0} + (x^{2})^{2}q + {\displaystyle \frac {16}{15}} 
(x^{2})^{5}q^{2} + {\displaystyle \frac {961}{420}} q^{3}(x
_{2})^{8} + {\displaystyle \frac {4105537}{623700}} q^{4}(x^{
2})^{11} + {\displaystyle \frac {291788599}{13097700}} q^{5}(
x^{2})^{14}+\cdots.
\label{hh}
\end{eqnarray}
If we invert the mirror maps and substitute them to (\ref{hh}), 
\begin{eqnarray}
&&w({\cal O}_{h}{\cal O}_{h}|(x^0(t^0,t^1,t^2),x^1(t^0,t^1,t^2),x^2(t^0,t^1,t^2)))_{0}=\no\\
&&=t^{0} + {\displaystyle \frac {1}{2}} (t^{2})^{2}Q + 
{\displaystyle \frac {1}{30}} (t^{2})^{5}Q^{2} + 
{\displaystyle \frac {3}{1120}} (t^{2})^{8}Q^{3} + 
{\displaystyle \frac {31}{124740}} (t^{2})^{11}Q^{4} + 
{\displaystyle \frac {1559}{62270208}} (t^{2})^{14}Q^{5}+\cdots\no\\
&&=t^{0}+\frac{1}{2!}(t^{2})^{2}Q+\frac{2^2}{5!}(t^{2})^{5}Q^{2}+\frac{3^2\cdot 12}{8!}(t^{2})^{8}Q^{3}
+\frac{4^2\cdot 620}{11!}(t^{2})^{11}Q^{4}+\frac{5^2\cdot 87304}{14!}(t^{2})^{14}Q^{5}+\cdots\no\\
&&\hspace{12cm}(Q:=e^{t^1}),
\end{eqnarray}
the result coincides with $\langle{\cal O}_{h}{\cal O}_{h}(t^0,t^1,t^2)\rangle_{0}$ 
computed from the associativity equation \cite{km}.
If we compute,
\begin{eqnarray}
&&w({\cal O}_{h^2}{\cal O}_{h^2}|(x^0,x^1,x^2))_{0}=\no\\
&&= q + {\displaystyle \frac {2}{3}} (x^{2})^{3}q^{2} + 
{\displaystyle \frac {17}{15}} q^{3}(x^{2})^{6} + 
{\displaystyle \frac {6455}{2268}} q^{4}(x^{2})^{9} + 
{\displaystyle \frac {4124497}{467775}} q^{5}(x^{2})^{12}+\cdots,
\end{eqnarray}
we obtain $\langle{\cal O}_{h^2}{\cal O}_{h^2}(t^0,t^1,t^2)\rangle_{0}$.
\begin{eqnarray}
&&w({\cal O}_{h^2}{\cal O}_{h^2}|(x^0(t^0,t^1,t^2),x^1(t^0,t^1,t^2),x^2(t^0,t^1,t^2)))_{0}=\no\\
&&=Q + {\displaystyle \frac {1}{6}} (t^{2})^{3}Q^{2} + 
{\displaystyle \frac {1}{60}} Q^{3}(t^{2})^{6} + 
{\displaystyle \frac {31}{18144}} Q^{4}(t^{2})^{9} + 
{\displaystyle \frac {1559}{8553600}} Q^{5}(t^{2})^{12}+\cdots\no\\
&&=Q+\frac{1}{3!}(t^{2})^{3}Q^{2}+\frac{12}{6!}(t^{2})^{6}Q^{3}
+\frac{620}{9!}(t^{2})^{9}Q^{4}+\frac{87304}{12!}(t^{2})^{12}Q^{5}+\cdots.
\end{eqnarray}

\subsection{Comparison with the Result of Iritani's $I$-function}
In this section, we demonstrate standard type of mirror computation of the $CP^2$-model following Iritani's work \cite{iri}.
We think that his method is fundamentally the same as the mirror computation by Barannikov \cite{bara}. According to 
Iritani \cite{irip}, this method starts from the following extended $I$-function.  
\begin{equation}
I_{CP^2}(z,h,y^1,y^2):=\sum_{n,m\geq 0}^{\infty}\exp (\frac{y^{1}h}{z})\cdot\frac{\displaystyle{\prod_{j=-\infty}^{0}(h+jz)^3\cdot\prod_{j=-\infty}^{0}(jz)}}
{\displaystyle{\prod_{j=-\infty}^{n-m}(h+jz)^2\cdot\prod_{j=-\infty}^{n}(h+jz)\cdot\prod_{j=-\infty}^{m}(jz)}}\cdot e^{n y^1}(y^2)^m.
\end{equation}
Here, $h$ is the hyperplane class of $CP^2$ and $h^3=0$. 
In this $I$-function, the parameter $z$ plays the role of $\hbar$ in Barannikov's formalism.
For a function $f(z,h,y_1,y_2)$ that includes only positive power of $h$, we denote by $(f)_{h^i}\;(i=0,1,2)$ the 
coefficient of $h^i$. We define here the following $3\times 3$matrix.    
\begin{equation}
S(z,y^1,y^2):=\left(\begin{array}{ccc}\bigl(I_{CP^2}\bigr)_{h^0}&\bigl(z\frac{\d}{\d y_1}I_{CP^2}\bigr)_{h^0} &
\bigl(z^2\frac{{\d}^2}{\d (y^1)^2}I_{CP^2}\bigr)_{h^0}\\
\bigl(I_{CP^2}\bigr)_{h^1}&\bigl(z\frac{\d}{\d y_1}I_{CP^2}\bigr)_{h^1} &
\bigl(z^2\frac{{\d}^2}{\d (y^1)^2}I_{CP^2}\bigr)_{h^1}\\
\bigl(I_{CP^2}\bigr)_{h^2}&\bigl(z\frac{\d}{\d y_1}I_{CP^2}\bigr)_{h^2} &
\bigl(z^2\frac{{\d}^2}{\d (y^1)^2}I_{CP^2}\bigr)_{h^2}
\end{array}\right).
\label{S1}
\end{equation}
This matrix has the following structure:
\begin{eqnarray}
&&S(z,y^1,y^2):=\left(\begin{array}{ccc}1&0&0\\
\frac{y^1}{z} &1&0\\
\frac{(y^1)^2}{2z^2}&\frac{y^1}{z}&1\end{array}
\right)\cdot M(z,y^1,y^2),\no\\
&&M(z,y^1,y^2)=\sum_{n,m\geq 0}^{\infty}M_{nm}(z)\cdot e^{ny^1}(y^2)^m,\no\\
&&M_{00}(z)=\left(\begin{array}{ccc}1&0&0\\
0&1&0\\
0&0&1\end{array}
\right),
\label{stS}
\end{eqnarray} 
where $M(z,y^1,y^2)$ and $M_{nm}(z)$ are $3\times 3$ matrices. By computing $M(z,y^1,y^2)$ explicitly, we can observe that it includes both positive and 
negative powers of $z$. Following Iritani, we factorize $M(z,y^1,y^2)$ into the form $M(z,y^1,y^2)=M_{-}(z,y^1,y^2)M_{+}(z,y^1,y^2)$ where 
$M_{-}(z,y^1,y^2)$ (resp. $M_{+}(z,y^1,y^2)$) includes only negative (resp. non-negative) powers of $z$. Since $M(z,y^1,y^2)$ has the structure of 
power series in $e^{y^1}$ and $y^2$ with top term an identity matrix, we can execute this operation systematically. 
We then introduce the matrix:
\begin{equation}
S_{-}(z,y^1,y^2):=\left(\begin{array}{ccc}1&0&0\\
\frac{y^1}{z} &1&0\\
\frac{(y^1)^2}{2z^2}&\frac{y^1}{z}&1\end{array}
\right)\cdot M_{-}(z,y^1,y^2).
\label{s-}
\end{equation}
The B-model connection matrix $C_{1}(y^1,y^2)$ of this setting is given by,
\begin{equation}
C_{1}(y^1,y^2)=\bigl(S_{-}(z,y^1,y^2)\bigr)^{-1}\cdot z\frac{\d}{\d y^1}S_{-}(z,y^1,y^2).
\end{equation}  
From the general theory of Iritani's method, $C_{1}(y^1,y^2)$ is free of $z$-dependence. Let $(A)_{ij} \;(i,j=1,2,3)$ be the $(i,j)$-element of $3\times 3$ matrix $A$.
Partial derivatives of the flat coordinates $t^i\;(i=0,1,2)$ associated with $h^i$ by the B-model coordinate $y^{1}$ are read off from $(C_{1}(y^1,y^2))_{i+1,1}$.
\begin{eqnarray}
&&\frac{\d t^0}{\d y^1}=(C_{1}(y^1,y^2))_{11}=\no\\
&&-{\displaystyle \frac {1}{2}}\tilde{q}(y^2)^{2}
+{\displaystyle \frac {13}{15}}\tilde{q}^{2}(y_2)^{5}
-{\displaystyle \frac {3167}{840}}\tilde{q}^{3}(y^2)^{8}
+{\displaystyle \frac {44552}{2079}}\tilde{q}^{4}(y^2)^{11}
-{\displaystyle \frac {450037373}{3243240}}\tilde{q}^{5}(y^2)^{14}+\cdots,\no\\
&&\frac{\d t^1}{\d y^1}=(C_{1}(y^1,y^2))_{21}=\no\\
&&1
+{\displaystyle \frac {1}{6}}\tilde{q}(y^2)^{3}
-{\displaystyle \frac {11}{15}}\tilde{q}^{2}(y^2)^{6} 
+{\displaystyle \frac {229}{56}}\tilde{q}^{3}(y^2)^{9}
-{\displaystyle \frac {775267}{29700}}\tilde{q}^{4}(y^2)^{12} 
+{\displaystyle \frac {233170937}{1289925}}\tilde{q}^{5}(y^2)^{15}+\cdots,\no\\ 
&&\frac{\d t^2}{\d y^1}=(C_{1}(y^1,y^2))_{31}=\no\\
&&-{\displaystyle \frac {5}{12}}\tilde{q}(y^2)^{4}
+{\displaystyle \frac {1241}{630}}\tilde{q}^{2}(y^2)^{7} 
-{\displaystyle \frac {47977}{4200}}\tilde{q}^{3}(y^2)^{10}
+{\displaystyle \frac{201402797}{2702700}}\tilde{q}^{4}(y^2)^{13}\,,
-{\displaystyle \frac {475054027589}{908107200}}\tilde{q}^{5}(y^2)^{16}+\cdots,\no\\
\label{irijac}
\end{eqnarray} 
where $\tilde{q}=e^{y^1}$.
By integrating the above equations in $y^1$, we obtain the mirror map,
\begin{eqnarray}
&&{t^{0}}= {y^{0}}
-{\displaystyle \frac {1}{2}}\tilde{q}(y^{2})^{2}
+{\displaystyle \frac {13}{30}}\tilde{q}^{2}(y^{2})^{5}
-{\displaystyle \frac {3167}{2520}}\tilde{q}^{3}(y^{2})^{8}
+{\displaystyle \frac {11138}{2079}}\tilde{q}^{4}(y^{2})^{11} 
-{\displaystyle \frac {450037373}{16216200}}\tilde{q}^{5}(y^{2})^{14} 
+\cdots\no\\
&&{t^{1}}= {y^{1}}
+{\displaystyle \frac {1}{6}}\tilde{q}(y^{2})^{3}
-{\displaystyle \frac {11}{30}}\tilde{q}^{2}(y^{2})^{6}
+{\displaystyle \frac {229}{168}}\tilde{q}^{3}(y^{2})^{9} 
-{\displaystyle \frac {775267}{118800}}\tilde{q}^{4}(y^{2})^{12} 
+{\displaystyle \frac {233170937}{6449625}}\tilde{q}^{5}(y^2)^{15}+\cdots\no\\ 
&&t^{2}= {y^2} 
-{\displaystyle \frac {5}{12}}\tilde{q}(y^{2})^{4}
+{\displaystyle \frac {1241}{1260}}\tilde{q}^{2}(y^2)^{7}
-{\displaystyle \frac {47977}{12600}}\tilde{q}^{3}(y^{2})^{10} 
+{\displaystyle \frac {201402797}{10810800}}\tilde{q}^{4}(y^{2})^{13} 
-{\displaystyle \frac {475054027589}{4540536000}}\tilde{q}^{5}(y^{2})^{16}+\cdots\no\\
\label{irimir}
\end{eqnarray}
On the other hand, the matrix elements $(C_{1}(y^1,y^2))_{12}=(C_{1}(y^1,y^2))_{23}, (C_{1}(y^1,y^2))_{13}, (C_{1}(y^1,y^2))_{22}$ give us information 
of the Gromov-Witten invariants. By integrating these matrix elements in $y^1$, we obtain the following functions.
\begin{eqnarray}
&&f_1=\tilde{q}{y^{2}} 
-{\displaystyle \frac {1}{6}}\tilde{q}^{2}(y^{2})^{4}
+{\displaystyle \frac {289}{630}}\tilde{q}^{3}(y^{2})^{7}
-{\displaystyle \frac {35873}{18900}}\tilde{q}^{4}(y^{2})^{10} 
+{\displaystyle \frac {156650191}{16216200}}\tilde{q}^{5}(y^{2})^{13}+\cdots,\no\\
&&f_2=\tilde{q}
+{\displaystyle \frac {1}{3}}\tilde{q}^{2}(y^{2})^{3}
-{\displaystyle \frac {22}{45}}\tilde{q}^{3}(y^{2})^{6}
+{\displaystyle \frac {1261}{756}}\tilde{q}^{4}(y^{2})^{9}
-{\displaystyle \frac {2405639}{311850}}\tilde{q}^{5}(y^{2})^{12}+\cdots,\no\\  
&&f_3=y^{0}
+{\displaystyle \frac {2}{15}}\tilde{q}^{2}(y^2)^{5}
-{\displaystyle \frac {613}{1260}}\tilde{q}^{3}(y^{2})^{8} 
+{\displaystyle \frac {4751}{2079}}\tilde{q}^{4}(y^{2})^{11} 
-{\displaystyle \frac {101313427}{8108100}}\tilde{q}^{5}(y^{2})^{14}+\cdots. 
\end{eqnarray}
If we expand these functions in $t^0, t^2$ and $Q=e^{t^1}$ by substituting the inversion of the mirror map, the result 
turns out to be, 
\begin{eqnarray}
&&f_1=Q{t^{2}}+{\displaystyle \frac {1}{12}}Q^2(t^{2})^{4}
+{\displaystyle \frac {1}{140}}Q^3(t^{2})^{7} + 
{\displaystyle \frac {31}{45360}}Q^4(t^{2})^{10}+ 
{\displaystyle \frac {1559}{22239360}}Q^5(t^{2})^{13}+\cdots,\no\\
&&f_2=Q+{\displaystyle \frac {1}{6}}Q^2(t^{2})^{3} + 
{\displaystyle \frac {1}{60}}Q^3(t^{2})^{6} + 
{\displaystyle \frac {31}{18144}}Q^4(t^{2})^{9}+ 
{\displaystyle \frac {1559}{8553600}}Q^5(t^2)^{12}+\cdots,\no\\
&&f_3=t^{0} + {\displaystyle \frac {1}{2}}Q(t^{2})^{2}+ 
{\displaystyle \frac {1}{30}}Q^2(t^{2})^{5}+ 
{\displaystyle \frac {3}{1120}}Q^3(t^{2})^{8}+ 
{\displaystyle \frac {31}{124740}}Q^4(t^{2})^{11}+ 
{\displaystyle \frac {1559}{62270208}}Q^5(t^{2})^{14}+\cdots.
\end{eqnarray}
Hence they reproduce $\langle{\cal O}_{h}{\cal O}_{h^2}(t^0,t^1,t^2)\rangle$, $\langle{\cal O}_{h^2}{\cal O}_{h^2}(t^0,t^1,t^2)\rangle$ and 
$\langle{\cal O}_{h}{\cal O}_{h}(t^0,t^1,t^2)\rangle$ respectively. The final results coincide with our computation, 
but we can see from (\ref{cp2mirr}) and (\ref{irimir}) that the mirror map in this case is different from our mirror map. 
As we have mentioned in Section 1, there exists infinitely many ways to include the parameter $y^2$ into the $I$-function \cite{irip}. 
Since we have tested only one possibility here, we cannot conclude that our formalism and Iritani's formalism have no connection.  
   
\section{Open String Case}
In this section, we discuss generalization of the multi-point virtual structure constants to the open string case. 
First, we consider anti-holomolphic involution $\varphi:CP^2\rightarrow CP^2$ defined by $\varphi(X_1:X_2:X_3)=
(\bar{X_1}:\bar{X_2}:\bar{X_3})$.  The subset invariant 
under $\varphi$ is $RP^2$, which is a Lagrangian submanifold of $CP^2$. Next, we pick up a quasi map from 
$CP^1$ to $CP^2$ of degree $2d-1$,
\begin{equation}
q_{2d-1}(s:t):=[\sum_{j=0}^{2d-1}{\bf a}_{j}s^{j}t^{2d-1-j}],\;\;({\bf a}_j\in{\bf C}^3).
\end{equation}
We also introduce an involution $u:CP^1\rightarrow CP^1$ defined by $u(s:t)=(\bar{t}:\bar{s})$. 
With this set-up, we define a ${\bf Z}_2$-action on the quasi map given by $q\mapsto \varphi\circ q_{2d-1}\circ u$. 
Since $\varphi(q_{2d-1}(u(s:t)))=[\overline{\sum_{j=0}^{2d-1}{\bf a}_{j}{\bar{t}}^{j}{\bar{s}}^{2d-1-j}}]
=[\sum_{j=0}^{2d-1}\bar{\bf a}_{2d-1-j}s^{j}t^{2d-1-j}]$, this action induces an involution on the parameter 
space of quasi maps,
\begin{equation}
({\bf a}_0,{\bf a}_1,\cdots,{\bf a}_{2d-1})\rightarrow (\bar{\bf a}_{2d-1},\bar{\bf a}_{2d-2},\cdots,\bar{\bf a}_0).
\label{pinvo}
\end{equation}
Let us consider a quasi map invariant under the above involution,
\begin{equation}
q_{2d-1}(s:t)=[\sum_{j=0}^{d-1}({\bf a}_js^{j}t^{2d-1-j}+\bar{\bf a}_js^{2d-1-j}t^{j}].
\label{qinv}
\end{equation}
It maps the equator of $CP^1$ $(\{(e^{\sqrt{-1}\theta}:1)\;|\;\theta\in[0,2\pi)\;\})$, which is invariant under $u$, to 
$RP^2$ because,
\begin{equation}
q_{2d-1}(1:e^{\sqrt{-1}\theta})=[\sum_{j=0}^{d-1}({\bf a}_je^{\sqrt{-1}j\theta}+\bar{\bf a}_je^{\sqrt{-1}(2d-1-j)\theta}]
=[\sum_{j=0}^{d-1}({\bf a}_je^{\sqrt{-1}(j-\frac{2d-1}{2})\theta}+\bar{\bf a}_je^{\sqrt{-1}(\frac{2d-1}{2}-j)\theta}].
\label{equa}
\end{equation}
Therefore, a quasi map invariant under the involution (\ref{pinvo}) can be regarded as a quasi map from upper half 
disk of $CP^1$ to $CP^2$, which maps boundary of the disk to the Lagrangian submanifold $RP^2$.

Next, we consider $Mp_{0,2|2n}(3,2d-1)$. It was defined by dividing the set,
\begin{equation}
Up_{0,2|2n}(3,2d-1):=\{\bigl(({\bf a}_0,\cdots,{\bf a}_{2d-1}),(z_1,\cdots,z_{2n})\bigr)\;|
\;{\bf a}_i\in{\bf C}^{3},\;z_i\in{\bf C}^{\times},\;{\bf a}_0,{\bf a}_{2d-1}\neq {\bf 0}\;\},
\end{equation}
by the two ${\bf C}^{\times}$ actions given in (\ref{action}). Motivated by the previous discussion, we introduce an 
involution $v: Up_{0,2|2n}(3,2d-1)\rightarrow Up_{0,2|2n}(3,2d-1)$ as follows.
\begin{equation}
v(({\bf a}_1,\cdots,{\bf a}_{2d-1}),(z_1,\cdots,z_{2n}))=\bigl((\bar{\bf a}_{2d-1}, \bar{\bf a}_{2d-2}, 
\cdots,\bar{\bf a}_{0}),(\frac{1}{\bar{z}_{2n}},\frac{1}{\bar{z}_{2n-1}},  \cdots,\frac{1}{\bar{z}_{1}})\bigr).
\label{mpinvo}
\end{equation}
It is easy to check that $v$ is compatible with equivalence relation by the two ${\bf C}^{\times}$ actions.
Hence it induces an involution $vp: Mp_{0,2|2n}(3,2d-1)\rightarrow Mp_{0,2|2n}(3,2d-1)$. We can easily extend 
$vp$ to whole  $\widetilde{Mp}_{0,2|2n}(3,2d-1)$ by looking back at the construction in Section 2.1.
Let us denote the extended involution by $\widetilde{vp}$. With this set-up, we define the moduli space 
$\widetilde{Mp}_{D,1|n}(CP^2/RP^2,2d-1)$ as the invariant subset of $\widetilde{Mp}_{0,2|2n}(3,2d-1)$ under 
$\widetilde{vp}$. Roughly speaking, the degrees of freedom of this moduli space are described by,
\begin{eqnarray} 
&&\{\bigl(({\bf a}_0,\cdots,{\bf a}_{d-1}),(z_1,\cdots,z_{n})\bigr)\;|\;{\bf a}_i\in {\bf C}^3,\;0<|z_i|\leq 1,\;{\bf a}_0\neq{\bf 0}\}/({\bf R}_{>0}\times U(1)),
\no\\
&& re^{\sqrt{-1}\theta}\cdot\bigl(({\bf a}_0,\cdots,{\bf a}_{d-1}),(z_1,\cdots,z_{n})\bigr)=\bigl((re^{\sqrt{-1}\frac{2d-1}{2}\theta}{\bf a}_0,
,\cdots,re^{\sqrt{-1}\frac{1}{2}\theta}{\bf a}_{d-1}),(e^{\sqrt{-1}\theta}z_1,\cdots,e^{\sqrt{-1}\theta}z_{n})\bigr),\no\\
&&\vspace{10cm}(re^{\sqrt{-1}\theta}\in{\bf R}_{>0}\times U(1)),
\label{ru1} 
\end{eqnarray}
 and they are half of the ones of $\widetilde{Mp}_{0,2|2n}(3,2d-1)$.
Evaluation map $ev_i:\widetilde{Mp}_{D,1|n}(CP^2/RP^2,2d-1)\rightarrow CP^2$ at the $i$-th marked point $z_i$ is defined 
in the same way as the closed string case. Note here that our construction allows the marked points to lie on the boundary of the disk.   

At this stage, we can define open version of the multi-point virtual structure constant.
\begin{equation}
w({\cal O}_{h^a}|\prod_{i=1}^{n}({\cal O}_{h^{m_i}}))_{disk,2d-1}:=
\int_{\widetilde{Mp}_{D,1|n}(CP^2/RP^2,2d-1)}ev_{0}^{*}(h^{a})\cdot\prod_{i=1}^{n}ev_{i}^{*}(h^{m_i}).
\label{defopw}
\end{equation}
Now, we compute the above intersection number by localization technique. First, we introduce $U(1)$ action flow on $\widetilde{Mp}_{D,1|n}(CP^2/RP^2,2d-1)$
which is induced from the following $U(1)$ action flow on the bulk part,
\begin{equation}
e^{\sqrt{-1}t}\cdot[\bigl(({\bf a}_0,\cdots,{\bf a}_{d-1}),(z_1,\cdots,z_{n})\bigr)]:=[\bigl((e^{\sqrt{-1}\theta_0t}{\bf a}_0,\cdots,
e^{\sqrt{-1}\theta_{d-1}t}{\bf a}_{d-1}),(z_1,\cdots,z_{n})\bigr)],\;\;(t\in {\bf R}).
\end{equation} 
We can fix the ambiguity coming from ${\bf R}_{>0}\times U(1)$ by regarding ${\bf a}_0$ as a point $[{\bf a}_0]$ in $CP^2$. 
Then we obtain,
\begin{eqnarray}
&&[\bigl((e^{\sqrt{-1}\theta_0t}{\bf a}_0,\cdots,
e^{\sqrt{-1}\theta_{d-1}t}{\bf a}_{d-1}),(z_1,\cdots,z_{n})\bigr)]=\no\\
&&\bigl(([{\bf a}_0], e^{\sqrt{-1}(\theta_1-\frac{2d-3}{2d-1}\theta_{0})t}{\bf a}_1,                        
\cdots,e^{\sqrt{-1}(\theta_{d-1}-\frac{1}{2d-1}\theta_0)t}{\bf a}_{d-1}),
( e^{\sqrt{-1}(-\frac{2}{2d-1}\theta_0)t}z_1,\cdots,e^{\sqrt{-1}(-\frac{2}{2d-1}\theta_0)t}z_{n})\bigr).\no\\
\label{u1tri}
\end{eqnarray}
Therefore, a point in the bulk part is fixed under the $U(1)$ action only if ${\bf a}_1={\bf a}_2=\cdots={\bf a}_{d-1}={\bf 0}$ 
and $n=0$.
The strata of $\widetilde{Mp}_{D,1|n}(CP^2/RP^2,2d-1)$ are labeled by ordered partitions $(d_1,d_2,\cdots,d_l)\in OP_d$
and ordered decompositions $(\coprod_{i=0}^{l-1}A_i)\coprod(\coprod_{i=1}^{l}B_i)=\{1,2,\cdots,n\}$, 
The stratum labeled by $(d_1,d_2,\cdots,d_l)$ and $(\coprod_{i=0}^{l-1}A_i)\coprod(\coprod_{i=1}^{l}B_i)=\{1,2,\cdots,n\}$ is 
given as follows.
\begin{eqnarray}
&&\widetilde{Mp}_{0,2||A_{0}|}(3,0)\mathop{\times}_{CP^{2}}
Mp_{0,2||B_1|}(3,d_1)\mathop{\times}_{CP^{2}}\widetilde{Mp}_{0,2||A_{1}|}(3,0)\mathop{\times}_{CP^{2}}
Mp_{0,2||B_2|}(3,d_2)\mathop{\times}_{CP^{2}}
\cdots\no\\
&&
\cdots
\mathop{\times}_{CP^{2}}\widetilde{Mp}_{0,2||A_{l-2}|}(3,0)
\mathop{\times}_{CP^{2}}Mp_{0,2||B_{l-1}|}(3,d_{l-1})\mathop{\times}_{CP^{2}}\widetilde{Mp}_{0,2||A_{l-1}|}(3,0)
\mathop{\times}_{CP^{2}}Mp_{D,1||B_l|}(CP^2/RP^2,2d_l-1) .\no\\
\label{opstr}
\end{eqnarray}
From the observation above and the discussion in Section 2.2, we can see that non empty fixed point set comes from the 
stratum that satisfy $B_i=\emptyset\;(i=1,2,\cdots,l)$. The fixed point set coming from the stratum 
labeled by  $(d_1,d_2,\cdots,d_l)$ and $\coprod_{i=0}^{l-1}A_i=\{1,2,\cdots,n\}$ is given by 
$\prod_{i=0}^{l-1}(CP^2)_i\times\prod_{i=0}^{l-1}\overline{M}_{0,2||A_i|}$ in the same way as the closed string case.
Now, we pay attention to the fixed point set of $Mp_{D,2|0}(CP^2/RP^2,2d-1)$ under the $U(1)$ action.
It is given by $\{([{\bf a}_0],{\bf 0},\cdots,{\bf 0})\}=CP^2$, but we have residual $U(1)$ action coming from ${\bf R}_{>0}\times U(1)$
, which is generated by $\zeta$ satisfying $\zeta^{\frac{2d-1}{2}}=1$. We formally interpret this action as the one caused by 
${\bf Z}_{\frac{2d-1}{2}}$, a quasi cyclic group of order $\frac{2d-1}{2}$. We proceed computation by assuming this quasi group.
Then we can regard the above $CP^2$ as the set of orbifold singularities on which ${\bf Z}_{\frac{2d-1}{2}}$ acts. With this 
consideration, we can also regard the fixed point set $\prod_{i=0}^{l-1}(CP^2)_i\times\prod_{i=0}^{l-1}\overline{M}_{0,2||A_i|}$
as the set of orbifold singularities on which $\bigl(\prod_{i=1}^{l-1}{\bf Z}_{d_i}\bigr)\times{\bf Z}_{\frac{2d_l-1}{2}}$ acts.
 
We turn into determination of the localized integrand to compute $w({\cal O}_{h^a}|\prod_{j=1}^{n}({\cal O}_{h^{m_j}}))_{disk,2d-1}$, that comes from the fixed point set labeled by $(d_1,d_2,\cdots,d_l)$ and $\coprod_{i=0}^{l-1}A_i=\{1,2,\cdots,n\}$.
But there are many overlaps with the discussion in the closed string case. Hence we determine the localized integrand only 
from the fixed point set labeled by $(d)$ and $A_0=\{1,2,\cdots,n\}$. In this case, the fixed point set is given by $(CP^2)_0\times \overline{M}_{0,2|n}$.

The normal bundle for this set comes from the following contributions.
\begin{itemize}
\item[(i)] deforming ${\bf a}_i$ $(i=1,2,\cdots,d-1)$ from ${\bf 0}$.
\item[(ii)] resolving the nodal singularity.  
\end{itemize} 
 
From the $U(1)$ action given in (\ref{ru1}), the part coming from (i) is identified with $\oplus_{j=1}^{d-1}({\cal O}_{CP^2}(\frac{2d-1-2j}{2d-1})^{\oplus 3})$. 
The $U(1)$ character for ${\cal O}_{CP^2}(\frac{2d-1-2j}{2d-1})$ can be read off from (\ref{u1tri}) and it is given by,
\begin{equation}
\frac{2d-1-2j}{2d-1}\sqrt{-1}\theta_0-\sqrt{-1}\theta_j.
\end{equation}  
Therefore, the contribution coming from (i) is given as follows.
\begin{equation}
\frac{1}{\prod_{j=1}^{d-1}\bigl(\frac{2d-1-2j}{2d-1}(h_0+\sqrt{-1}\theta_0)-\sqrt{-1}\theta_j\bigr)^3},
\end{equation} 
where $h_0$ is the hyperplane class of $(CP^2)_0$. The part coming from (ii) is identified with $\frac{d}{d(\frac{s}{t})}\otimes T_{\infty}^{\prime}C_0$. 
Here, $C_0$ is the genus $0$ stable curve described by $\overline{M}_{0,2|n}$.
$\frac{d}{d(\frac{s}{t})}$ is identified with ${\cal O}_{CP^2}(\frac{2}{2d-1})$ and its $U(1)$ character is given by $\frac{2}{2d-1}\sqrt{-1}\theta_0$.
The contribution coming from this part turns out to be, 
\begin{equation}
\frac{1}{\frac{2}{2d-1}(h_0+\sqrt{-1}\theta_0 )+c_1(T_{\infty}^{\prime}C_0)}. 
\end{equation} 
As in the closed string case, we integrate the above equivariant form on $\overline{M}_{0,2|n}$. The result is, 
\begin{equation}
\biggl(\frac{2d-1}{2(h_0+\sqrt{-1}\theta_0 )}\biggr)^{n}.
\end{equation} 
The contribution from $ev_i^{*}(h^{m_i})$ (resp. $ev_0^{*}(h^a)$) is given by $(h_0+\sqrt{-1}\theta_0 )^{m_i}$ (resp. $(h_0+\sqrt{-1}\theta_0 )^{a}$).
With this set-up, we cam write down the localized integrand to compute $w({\cal O}_{h^a}|\prod_{i=1}^{n}({\cal O}_{h^{m_i}}))_{disk,2d-1}$,
that comes from the fixed point set labeled by $(2d-1)$ and $A_0=\{1,2,\cdots,n\}$.
\begin{equation}
\frac{2}{2d-1}\cdot\frac{1}{\prod_{j=1}^{d-1}\bigl(\frac{2d-1-2j}{2d-1}(h_0+\sqrt{-1}\theta_0)-\sqrt{-1}\theta_j\bigr)^3}\cdot(h_0+\sqrt{-1}\theta_0 )^{a}
\prod_{i=1}^{n}\biggl(\frac{2d-1}{2}(h_0+\sqrt{-1}\theta_0 )^{m_i-1}\biggr),
\end{equation}
where the factor $\frac{2}{2d-1}$ at the left end comes from the ${\bf Z}_{\frac{2d-1}{2}}$ action.
Remaining computation goes in the same way as the closed string case. The result of localization computation is given as follows.
\begin{eqnarray}
&&w({\cal O}_{h^a}|\prod_{i=1}^{n}{\cal O}_{h^{m_i}})_{disk,2d-1}=\sum_{(d_1,\cdots,d_l)\in OP_d}\frac{2}{2d_{l}-1}\biggl(\prod_{j=1}^{l-1}\frac{1}{d_j}\biggr)
\int_{(CP^{N-1})_0}\int_{(CP^{N-1})_1}\cdots\int_{(CP^{N-1})_{l-1}}\times\no\\
&&(h_{0}+\sqrt{-1}\theta_0)^a\prod_{j=1}^{l-1}\frac{1}{\displaystyle{\prod_{i=1}^{d_j-1}\biggl
(\bigl(\frac{d_{j}-i}{d_{j}}\bigr)(h_{f_{j-1}}+\sqrt{-1}\theta_{f_{j-1}})+\bigl( \frac{i}{d_{j}}\bigr)(h_{f_{j}}+\sqrt{-1}\theta_{f_{j}})-
\sqrt{-1}\theta_{f_{j-1}+i}
\biggr)^{3}}}\times\no\\
&&\prod_{j=1}^{l-2}\frac{1}{\displaystyle{\frac{h_{f_j}+\sqrt{-1}\theta_{f_j}-h_{f_{j-1}}-\sqrt{-1}\theta_{f_{j-1}}}{d_{j}}+
\frac{h_{f_j}+\sqrt{-1}\theta_{f_j}-h_{f_{j+1}}-\sqrt{-1}\theta_{f_{j+1}}}{d_{j+1}}}}\times\\
&&\frac{1}{\displaystyle{\prod_{i=1}^{d_l-1}\biggl
(\bigl(\frac{2d_{l}-1-2i}{2d_{l}-1}\bigr)(h_{f_{l-1}}+\sqrt{-1}\theta_{f_{l-1}})-
\sqrt{-1}\theta_{f_{l-1}+i}
\biggr)^{3}}}\times\no\\
&&
\frac{1}{\displaystyle{\frac{h_{f_{l-1}}+\sqrt{-1}\theta_{f_{l-1}}-h_{f_{l-2}}-\sqrt{-1}\theta_{f_{l-2}}}{d_{l-1}}+
\frac{2(h_{f_{l-1}}+\sqrt{-1}\theta_{f_{l-1}})}{2d_{l}-1}}}\times\no\\
&&\prod_{i=1}^{n}\bigl(\sum_{j=1}^{l-1}w_{m_i}^{d_j}(h_{f_{j-1}}+\sqrt{-1}\theta_{f_{j-1}},h_{f_j}+\sqrt{-1}\theta_{f_j})
+\frac{2d_l-1}{2}(h_{l-1}+\sqrt{-1}\theta_{l-1} )^{m_i-1}\bigr),
\label{wints}
\end{eqnarray} 
where we also used the alternate notation $0=f_0<f_1<\cdots<f_{l-1}<f_l=d$ for the ordered partition $(d_1,d_2,\cdots,d_l)$. 
We then use the trick of residue integral and non-equivariant limit $\theta_{j}\rightarrow 0$ $(j=0,1,\cdots,d-1)$ and obtain the 
following theorem.
\begin{theorem}
\begin{eqnarray}
w({\cal O}_{h^a}|\prod_{i=1}^{n}{\cal O}_{h^{m_i}})_{disk,2d-1}&=&
\cdot\frac{1}{(2\pi\sqrt{-1})^{d}}\oint_{E^0_{(0)}}\frac{dz_0}{(z_0)^3}\oint_{E^{1}_{(0)}}\frac{dz_1}{(z_1)^3}\cdots \oint_{E^{d-1}_{(0)}}\frac{dz_{d-1}}{(z_{d-1})^3}
\times\no\\
&&2(z_{0})^{a}\biggl(\prod_{j=1}^{d-1}\frac{1}{(2z_j-z_{j-1}-z_{j+1})}\biggr)\cdot
\prod_{i=1}^{n}\biggl(\sum_{j=1}^{d-1}w_{m_i}^1(z_{j-1},z_j)+\frac{1}{2}(z_{d-1})^{m_i-1}\biggr), \no\\
&&\hspace{8cm}(d\geq1),
\label{thm2}
\end{eqnarray}
where we formally set $z_{d}:=-z_{d-1}$. $\oint_{E^0_{(0)}}dz_0$ is the operation of taking a residue at $z_{0}=0$, $\oint_{E^{j}_{(0)}}dz_j$ 
$(j=1,\cdots d-2)$ is the operation of taking residues at $z_j=0, \frac{z_{j-1}+z_{j+1}}{2}$ and $\oint_{E^{d-1}_{(0)}}dz_{d-1}$ 
is the operation of taking residues at $z_{d-1}=0, \frac{z_{d-2}}{3}$ respectively.
\end{theorem}
Now, we can compute the open multi-point virtual structure constant $w({\cal O}_{h^a}|\prod_{i=1}^{n}{\cal O}_{h^{m_i}})_{disk,2d-1}$
by using the above formula. Next step is to answer the question whether or not we can compute open Gromov-Witten invariant
$\langle\prod_{i=1}^{n}{\cal O}_{h^{m_i}}\rangle_{disk,2d-1}$ of $CP^2$ from the virtual structure constants. 
For this purpose, we prepare some numerical data of open Gromov-Witten invariants of $CP^2$. In \cite{opv}, we 
proposed formulas to compute open Gromov-Witten invariants of degree $k$ hypersurface of $CP^{N-1}$ $(k:odd)$.
Here, we write down the formulas again.
\begin{prop}\cite{opv}
 The A-model amplitude $\langle\prod_{i=1}^{n}{\cal O}_{h^{a_i}}\rangle_{disk,2d-1}$ up to the $d=3$ case 
is given by sum of the following residue integrals.   
\begin{eqnarray}
\langle\prod_{i=1}^{n}{\cal O}_{h^{m_i}}\rangle_{disk,1}&=&\frac{1}{(2\pi\sqrt{-1})}\oint_{C_{0}}\frac{dz_0}{(z_{0})^{N}}
f^{N,k}_{1}(z_{0})\cdot2z_0\cdot\prod_{i=1}^{n}(\frac{(z_0)^{m_i-1}}{2}),\no\\
\langle\prod_{i=1}^{n}{\cal O}_{h^{m_i}}\rangle_{disk,3}&=&\frac{1}{(2\pi\sqrt{-1})}\oint_{C_{0}}\frac{dz_0}{(z_{0})^{N}}
f^{N,k}_{3}(z_{0})\cdot\frac{2z_0}{3}\cdot\prod_{i=1}^{n}(\frac{3(z_0)^{m_i-1}}{2})\no\\
&+&\frac{1}{(2\pi\sqrt{-1})^2}\oint_{C_{0}}\frac{dz_0}{(z_{0})^{N}}\oint_{C_{1}}\frac{dz_1}{(z_{1})^{N}}
f^{N,k}_{1}(z_{0})e^{k}(z_0,z_1)\frac{z_1-z_0}{kz_0(3z_0-z_1)}
\cdot\prod_{i=1}^{n}(\frac{(z_0)^{m_i-1}}{2}+w_{m_i}^1(z_0,z_1)), \no\\
\langle\prod_{i=1}^{n}{\cal O}_{h^{m_i}}\rangle_{disk,5}&=&\frac{1}{(2\pi\sqrt{-1})}\oint_{C_{0}}\frac{dz_0}{(z_{0})^{N}}
f^{N,k}_{5}(z_{0})\cdot\frac{2z_0}{5}\cdot\prod_{i=1}^{n}(\frac{5(z_0)^{m_i-1}}{2})\no\\
&+&\frac{1}{(2\pi\sqrt{-1})^2}\oint_{C_{0}}\frac{dz_0}{(z_{0})^{N}}\oint_{C_{1}}\frac{dz_1}{(z_{1})^{N}}
f^{N,k}_{3}(z_{0})e^{k}(z_0,z_1)\frac{z_1-z_0}{kz_0(\frac{5}{3}z_0-z_1)}
\cdot \prod_{i=1}^{n}(\frac{3(z_0)^{m_i-1}}{2}+w_{m_i}^1(z_0,z_1))\no\\
&+&\frac{1}{(2\pi\sqrt{-1})^3}\oint_{C_{0}}\frac{dz_0}{(z_{0})^{N}}\oint_{C_{1}}\frac{dz_1}{(z_{1})^{N}}
\oint_{C_{2}}\frac{dz_2}{(z_{2})^{N}}
f^{N,k}_{1}(z_{0})e^{k}(z_0,z_1)e^{k}(z_1,z_2)\times\no\\
&&\frac{z_2-z_1}{kz_0(3z_0-z_1)kz_1(2z_1-z_0-z_2)}
\prod_{i=1}^{n}(\frac{(z_0)^{m_i-1}}{2}+w_{m_i}^1(z_0,z_1)+w_{m_i}^1(z_1,z_2))\no\\
&+&\frac{1}{2}\frac{1}{(2\pi\sqrt{-1})^3}\oint_{C_{0}}\frac{dz_0}{(z_{0})^{N}}\oint_{C_{1}}\frac{dz_1}{(z_{1})^{N}}
\oint_{C_{2}}\frac{dz_2}{(z_{2})^{N}}
f^{N,k}_{1}(z_{0})e^{k}(z_0,z_1)e^{k}(z_0,z_2)\frac{1}{(kz_0)^2(2z_0)}\times\no\\
&&\prod_{i=1}^{n}(\frac{(z_0)^{m_i-1}}{2}+w_{m_i}^1(z_0,z_1)+w_{m_i}^1(z_0,z_2)),
\label{ogres}
\end{eqnarray}
where
\begin{eqnarray}
e^{k}(z,w):=\prod_{j=0}^{k}(jz+(k-j)w). \label{e^k},
\label{ek}
\end{eqnarray}
and, 
\begin{eqnarray}
f^{N,k}_{2d-1}(z)&:=&\frac{2}{2d-1}\cdot\frac{\displaystyle{\prod_{j=0}^{kd-\frac{k+1}{2}}(\frac{j(-z)+(k(2d-1)-j)z}{2d-1})}}
{\displaystyle{\prod_{j=1}^{d-1}(\frac{j(-z)+(2d-1-j)z}{2d-1})^N}}.
\label{def1}
\end{eqnarray}
In the above formulas, we take the residue integrals in ascending order of the subscript $i$ of $z_i$. $\frac{1}{2\sqrt{-1}}\oint_{C_i}dz_i$
means that we take the residues at $z_{i}=0,\frac{z_{i-1}+z_{i+1}}{2}$ (resp. $z_i=0$) if the integrand contains the factor
$\frac{1}{2z_i-z_{i-1}-z_{i+1}}$ (resp. otherwise).
\end{prop}
This proposition followed from the localization computation applied to the open Gromov-Witten invariants \cite{kont, walcher} and the non-equivariant limit.
If we set $N=4, k=1$, we can compute the open Gromov Witten invariants of $CP^2$. We then obtain the following data.
\begin{equation}
\langle{\cal O}_{h^2}\rangle_{disk,1}=2,\;\langle({\cal O}_{h^2})^4\rangle_{disk,3}=-\frac{9}{4}, \;\;\langle({\cal O}_{h^2})^7\rangle_{disk,5}=\frac{3361}{32}.
\label{opn}
\end{equation}
Our first approach to the question is to consider the generating function,
\begin{eqnarray}
w({\cal O}_{h^a}|(x^0,x^1,x^2))_{disk}&:=&
\sum_{d\geq 1,m_j\geq 0}w({\cal O}_{h^a}|\prod_{j=0}^{2}({\cal O}_{h^{j}})^{m_j})_{disk,2d-1}\cdot\prod_{j=0}^{\infty}\frac{(x^{j})^{m_j}}{m_j!},
\end{eqnarray} 
and to compute $w({\cal O}_{h^a}|(x^0(t^0,t^1,t^2),x^1(t^0,t^1,t^2),x^2(t^0,t^1,t^2)))_{disk}$ by using the mirror map (\ref{cp2mirr}).
But this naive approach did not reproduce the above data. With some trials and errors, we found that the r.h.s. of (\ref{thm2}) 
produces non-zero rational numbers when we formally insert ${\cal O}_{m_i}$ with $m_i\geq 3$. This fact led us to a 
new approach to consider insertions of ${\cal O}_{h^j}$ $(j=0,1,2,3,4,\cdots)$ even though $H^{*,*}(CP^2)$ is spanned by $1,h,h^2$.
Explicitly, we consider the generating function,
\begin{eqnarray}
w({\cal O}_{h^a}|(x^0,x^1,x^2,x_3,\cdots))_{disk}&:=&
\sum_{d\geq 1,m_j \geq 0}w({\cal O}_{h^a}|\prod_{j=0}^{\infty}({\cal O}_{h^{j}})^{m_j})_{disk,2d-1}\cdot\prod_{j=0}^{\infty}\frac{(x^{j})^{m_j}}{m_j!},
\end{eqnarray}
where $w({\cal O}_{h^a}|\prod_{j=0}^{\infty}({\cal O}_{h^{j}})^{m_j})_{disk,2d-1}$ is defined by, 
\begin{eqnarray}
w({\cal O}_{h^a}|\prod_{j=0}^{\infty}({\cal O}_{h^{j}})^{m_j})_{disk,2d-1}&=&
\cdot\frac{1}{(2\pi\sqrt{-1})^{d}}\oint_{E^0_{(0)}}\frac{dz_0}{(z_0)^3}\oint_{E^{1}_{(0)}}\frac{dz_1}{(z_1)^3}\cdots \oint_{E^{d-1}_{(0)}}\frac{dz_{d-1}}{(z_{d-1})^3}
\times\no\\
&&2(z_{0})^{a}\biggl(\prod_{j=1}^{d-1}\frac{1}{(2z_j-z_{j-1}-z_{j+1})}\biggr)\cdot
\prod_{j=0}^{\infty}\biggl(\sum_{i=1}^{d-1}w_{j}^1(z_{i-1},z_i)+\frac{1}{2}(z_{d-1})^{j-1}\biggr)^{m_j}. \no\\
\label{opd2}
\end{eqnarray}
At this stage, we have to define the mirror map for the variables $x^{j}$ $(j=0,1,2,3,\cdots)$. If we pay attention to the fact that the mirror map in 
the closed string case was given by, 
\begin{equation}
t^j(x_0,x_1,x_2)=w({\cal O}_{h^{2-j}}{\cal O}_{1}|(x^0,x^1,x^2))_{0}, 
\end{equation}
we are naturally led to the following definition,
\begin{eqnarray}
&&t^{j}(x^{*})=t^{j}(x^0,x^1,x^2,\cdots):=w({\cal O}_{h^{2-j}}{\cal O}_{1}|(x^0,x^1,x^2,\cdots))_{0}:=\no\\
&&x^j+
\sum_{d>0,m_j\geq 0}w({\cal O}_{h^{2-j}}{\cal O}_{1}|
\prod_{j=0}^{\infty}({\cal O}_{h^{j}})^{m_j})_{0,d}\cdot\prod_{j=0}^{\infty}\frac{(x^{j})^{m_j}}{m_j!}.
\end{eqnarray}
Here, $w({\cal O}_{h^{2-j}}{\cal O}_{1}|
\prod_{j=0}^{\infty}({\cal O}_{h^{j}})^{m_j})_{0,d}$ is defined by, 
\begin{eqnarray}
w({\cal O}_{h^{2-j}}{\cal O}_{1}|\prod_{j=0}^{\infty}({\cal O}_{h^{j}})^{m_j})_{0,d}&=&
\frac{1}{(2\pi\sqrt{-1})^{d+1}}\oint_{E^0_{(0)}}\frac{dz_0}{(z_0)^3}\oint_{E^1_{(0)}}\frac{dz_1}{(z_1)^3}\cdots \oint_{E^d_{(0)}}\frac{dz_d}{(z_d)^3}
\times\no\\
&&(z_{0})^{2-j}\cdot\biggl(\prod_{j=1}^{d-1}\frac{1}{(2z_j-z_{j-1}-z_{j+1})}\biggr)\cdot
\prod_{j=0}^{\infty}\biggl(\sum_{i=1}^{d}w_{j}^1(z_{i-1},z_i)\biggr)^{m_j}, \;\;(d>0).\no\\
\label{p2cmod}
\end{eqnarray}
This formula produces non-trivial rational number even when $2-j<0$!
With this set-up's, we conjecture that generating function of the open Gromov-Witten invariants of $CP^2$: 
\begin{eqnarray}
\langle{\cal O}_{h^a}(t^0,t^1,t^2,\cdots)\rangle_{disk}&:=&
\sum_{d\geq 1,\{m_j\}}\langle{\cal O}_{h^a}\prod_{j=0}^{\infty}({\cal O}_{h^{j}})^{m_j}\rangle_{disk,2d-1}\cdot\prod_{j=0}^{\infty}\frac{(t^{j})^{m_j}}{m_j!},
\label{opgene}
\end{eqnarray}
can be computed by the equality:
\begin{equation}
w({\cal O}_{h^a}|(x^0,x^1,x^2,\cdots))_{disk}=\langle{\cal O}_{h^a}(t^0(x^{*}),t^1(x^{*}),t^2(x^{*}),\cdots)\rangle_{disk},
\end{equation}
or, conversely, 
\begin{equation}
\langle{\cal O}_{h^a}(t^0,t^1,t^2,\cdots)\rangle_{disk}=w({\cal O}_{h^a}|(x^0(t^{*}),x^1(t^{*}),x^2(t^{*}),\cdots))_{disk}.
\label{opconj}
\end{equation}
We have one subtle remark here. Even when we fix $d$, the sum $\sum_{m_j \geq 0}w({\cal O}_{h^a}|\prod_{j=0}^{\infty}({\cal O}_{h^{j}})^{m_j})_{disk,2d-1}\cdot\prod_{j=0}^{\infty}\frac{(x^{j})^{m_j}}{m_j!}$ contains infinite terms 
because we have ${\cal O}_1$ and ${\cal O}_{h}$ insertions. As for ${\cal O}_{h}$ insertion, we have the equality,
\begin{equation}
w({\cal O}_{h^a}|{\cal O}_h\prod_{j=0}^{\infty}({\cal O}_{h^{j}})^{m_j})_{disk,2d-1}=(d-\frac{1}{2})
w({\cal O}_{h^a}|\prod_{j=0}^{\infty}({\cal O}_{h^{j}})^{m_j})_{disk,2d-1}.
\end{equation}
Therefore, we can simplify the above sum to,
\begin{equation}
\sum_{m_j \geq 0, (j\neq 1)}w({\cal O}_{h^a}|\prod_{j\geq 0, j\neq 1}({\cal O}_{h^{j}})^{m_j})_{disk,2d-1}e^{(d-\frac{1}{2})x^1}\prod_{j\geq0,j\neq1}
\frac{(x^{j})^{m_j}}{m_j!}.
\label{opsum}
\end{equation} 
But unlike the closed string case, ${\cal O}_1$ insertion does not kill the open virtual structure constants.  We have infinite summands 
in the sum (\ref{opsum}). Hence our conjecture includes infinite procedures of computation even when we compute the generating function (\ref{opgene})
up to some finite $d$. So we limit our intention to compute $\langle({\cal O}_{h^2})^{3d-2}\rangle_{disk,2d-1}$, which are only non-trivial open 
Gromov-Witten invariants for $CP^2$ in usual sense.  To obtain $\langle({\cal O}_{h^2})^{3d-2}\rangle_{disk,2d-1}$ up to fixed $d$, we found that 
we can truncate the parameters to $x^0,x^1,\cdots,x^{d}$ and the number of ${\cal O}_1$ insertions $m_0$ in 
$\sum_{m_j \geq 0, (j\neq 1)}w({\cal O}_{h^a}|\prod_{j\geq 0, j\neq 1}({\cal O}_{h^{j}})^{m_j})_{disk,2f-1}e^{(f-\frac{1}{2})x^1}\prod_{j\geq0,j\neq1}
\frac{(x^{j})^{m_j}}{m_j!}$ to $0,1,\cdots,d-f$. We observed that these truncations do not affect the result of computation given in 
(\ref{opconj}) up to degree $2d-1$ and $t^0,t^1,\cdots,t^{d}$. Of course, even after this truncation, the size of computation is huge. 
So, we only demonstrate the computation up to $d=3$ here. The mirror map for $t^0,t^1,t^2,t^3$ is given as follows.
   
\begin{eqnarray}
{t_{3}} &:=& (x^{3}) + q({\displaystyle \frac {1}{12}} (x^{2})
^{5} + {\displaystyle \frac {7}{6}} (x^{2})^{3}x^{3} + 
{\displaystyle \frac {5}{2}} (x^{3})^{2}x^{2})+ q^{2}({\displaystyle \frac {73}{336}} (x^{2})^{8}
 + {\displaystyle \frac {64}{15}} x^{3}(x^{2})^{6} + 
{\displaystyle \frac {181}{8}} (x^{3})^{2}(x^{2})^{4} + \no\\
&&+{\displaystyle \frac {97}{3}} (x^{2})^{2}(x^{3})^{3} + 
{\displaystyle \frac {35}{6}} (x^{3})^{4})+\cdots,\no \\
&&\no\\
{t_{2}} &:=& (x^{2}) + q({\displaystyle \frac {1}{4}} (x^{2})^{
4} + {\displaystyle \frac {3}{2}} (x^{3})^{2} + 2(x^{2})^{2}
x^{3}) + q^{2}({\displaystyle \frac {33}{70}} (x^{2})^{7}
 + {\displaystyle \frac {203}{30}} x^{3}(x^{2})^{5} + 
{\displaystyle \frac {47}{2}} (x^{2})^{3}(x^{3})^{2} + 17{x
_{3}}^{3}x^{2})+\cdots,\no\\
&&\no\\
{t_{1}} &:=& {x_{1}} + ({\displaystyle \frac {1}{2}} (x^{2})^{3}
 + 2x^{2}x^{3})q + ({\displaystyle \frac {7}{10}} {x
_{2}}^{6} + {\displaystyle \frac {22}{3}} x^{3}(x^{2})^{4}
 + {\displaystyle \frac {61}{4}} (x^{2})^{2}(x^{3})^{2} + 
{\displaystyle \frac {7}{2}} (x^{3})^{3})q^{2}+\cdots,\no\\
{t_{0}} &:=& (x^{0}) + ({\displaystyle \frac {1}{2}} (x^{2})^{2}
 + x^{3})q + ({\displaystyle \frac {8}{15}} (x^{2})^{5} + 
{\displaystyle \frac {13}{3}} (x^{2})^{3}x^{3} + 
{\displaystyle \frac {11}{2}} (x^{3})^{2}x^{2})q^{2}+\cdots.
\label{opex}
\end{eqnarray}
To obtain $\langle({\cal O}_{h^2})^{3d-2}\rangle_{disk,2d-1}$, it is enough to compute,
\begin{eqnarray}
&&w({\cal O}_{h}|(x^0,x^1,x^2,x^3))_{disk}=\no\\
&&= ((x^{2}) + {\displaystyle \frac {1}{2}} x^{0}
x^{3} + {\displaystyle \frac {1}{8}} x^{0}(x^{2})^{2}
 + {\displaystyle \frac {1}{16}} (x^{0})^{2}x^{2}x^{3}
 + {\displaystyle \frac {1}{192}} (x^{0})^{2}(x^{2})^{3} + 
{\displaystyle \frac {1}{192}} (x^{0})^{3}(x^{3})^{2} + 
{\displaystyle \frac {1}{384}} (x^{0})^{3}(x^{2})^{2}x^{3}+
 \no\\
&&+ {\displaystyle \frac {1}{9216}} (x^{0})^{3}(x^{2})^{4})q^{(1/2)}+
({\displaystyle \frac {5}{4}} (x^{3})^{2} + {\displaystyle 
\frac {21}{8}} (x^{2})^{2}x^{3} + {\displaystyle \frac {27
}{64}} (x^{2})^{4} + {\displaystyle \frac {99}{1280}} x^{0}
(x^{2})^{5} + {\displaystyle \frac {27}{32}} (x^{2})^{3}x^{3}x^{0} +\no\\ 
&&+{\displaystyle \frac {21}{16}} x^{2}
(x^{3})^{2}(x^{0}) + {\displaystyle \frac {117}{20480}} (x^{0})^{2}
(x^{2})^{6})q^{(3/2)} + ({\displaystyle \frac {246023}{322560}} 
(x^{2})^{7} + {\displaystyle \frac {7489}{768}} (x^{2})^{5}
x^{3} + {\displaystyle \frac {1889}{64}} (x^{2})^{3}(x^{3})^{2}
 + \no\\
&&+{\displaystyle \frac {833}{48}} x^{2}(x^{3})^{3})q^{(5/2)}+\cdots,\;\;(q=e^{x^1}),
\label{w1}
\end{eqnarray}
but we also computed,
\begin{eqnarray} 
&&w({\cal O}_{h^2}|(x^0,x^1,x^2,x^3))_{disk} =\no\\
&&= (2 + {\displaystyle \frac {1}{2}} x^{0}
x^{2} + {\displaystyle \frac {1}{8}} x^{3}(x^{0})^{2} + 
{\displaystyle \frac {1}{32}} (x^{0})^{2}(x^{2})^{2} + 
{\displaystyle \frac {1}{96}} (x^{0})^{3}x^{2}x^{3} + 
{\displaystyle \frac {1}{1152}} (x^{0})^{3}(x^{2})^{3})
q^{(1/2)}+ \no\\
&&+ ({\displaystyle \frac {3}{2}} x^{2}x^{3} + 
{\displaystyle \frac {3}{8}} (x^{2})^{3} + 
{\displaystyle \frac {9}{128}} (x^{2})^{4}x^{0} + 
{\displaystyle \frac {9}{16}} (x^{2})^{2}x^{3}x^{0} + 
{\displaystyle \frac {3}{8}} (x^{3})^{2}x^{0} + 
{\displaystyle \frac {27}{5120}} (x^{2})^{5}(x^{0})^{2})q^{(3/2)}+ \no\\
&&+ ({\displaystyle \frac {12823}{23040}} (x^{2})^{6} + 
{\displaystyle \frac {2219}{384}} (x^{2})^{4}x^{3} + 
{\displaystyle \frac {391}{32}} (x^{2})^{2}(x^{3})^{2} + 
{\displaystyle \frac {67}{24}} (x^{3})^{3})q^{(5/2)}+\cdots,
\label{w2}
\end{eqnarray}
to check integrable condition.
By inverting the mirror map (\ref{opex}) and substituting the result to (\ref{w1}) and (\ref{w2}), we obtain the generating function 
of open Gromov-Witten invariants,
\begin{eqnarray}
&& \langle{\cal O}_{h}(t^0,t^1,t^2,t^3)\rangle_{disk}|_{t^0=0}=w({\cal O}_{h}|(x^0(t^{*}),x^1(t^{*}),x^2(t^{*}),\cdots))_{disk}=\no\\ 
&&= Q^{(1/2)}(t^{2})+ ( - {\displaystyle \frac {3}{4}} (t^{3})^{2} - 
{\displaystyle \frac {3}{4}} (t^{2})^{2}t^{3} - 
{\displaystyle \frac {9}{64}} (t^{2})^{4})Q^{(3/2)}+({\displaystyle \frac {3361}{64512}} (t^{2})^{7}
 + {\displaystyle \frac {33}{64}} (t^{2})^{5}t^{3} + \no\\
 &&+{\displaystyle \frac {145}{96}} (t^{3})^{2}(t^{2})^{3} + 
{\displaystyle \frac {65}{48}} (t^{3})^{3}t^{2})Q^{(5/2)}+\cdots, \;\;(Q=e^{t^1}),
\end{eqnarray}
and, 
\begin{eqnarray}
&&\langle{\cal O}_{h^2}(t^0,t^1,t^2,t^3)\rangle_{disk}|_{t^0=0}=w({\cal O}_{h^2}|(x^0(t^{*}),x^1(t^{*}),x^2(t^{*}),\cdots))_{disk}=\no\\ 
&&=  2Q^{(1/2)}+ ( - {\displaystyle \frac {3}{8}} (t^{2})^{3} - t^{2}
t^{3})Q^{(3/2)} +({\displaystyle \frac {13}{24}} (t^{3})^{3} + 
{\displaystyle \frac {3361}{23040}} (t^{2})^{6} + 
{\displaystyle \frac {33}{32}} (t^{2})^{4}t^{3} + 
{\displaystyle \frac {29}{16}} (t^{3})^{2}(t^{2})^{2})Q^{(5/2)}+\cdots.\no\\
\end{eqnarray}
Here we set the variable $t^0$ to $0$ to simplify the formulas.
Note that the integrable condition,
\begin{equation}
\frac{\d}{\d t^2}\langle{\cal O}_{h}(t^0,t^1,t^2,t^3)\rangle_{disk}=\frac{\d}{\d t^1}\langle{\cal O}_{h^2}(t^0,t^1,t^2,t^3)\rangle_{disk},
\end{equation}
is satisfied. The numerical data (\ref{opn}) are also reproduced. We extended the computation up to $d=6$ and obtained the Table in Section 1.

\section{Application to General Type Projective Hypersurface} 
In this section, we discuss application of the multi-point virtual structure constants to Gromov-Witten
invariants of degree $k$ hypersurface in $CP^{N-1}$ (in our works, we denote it by $M_{N}^{k}$).
Let $w^{N,k}({\cal O}_{h^a}{\cal O}_{h^{b}}|\prod_{i=1}^{n}{\cal O}_{h^{m_i}})_{0,d}$ and $w^{N,k}({\cal O}_{h^a}|\prod_{i=1}^{n}{\cal O}_{h^{m_i}})_{disk,2d-1}$
be closed and open multi-point virtual structure constants for $M_{N}^{k}$. They are defined as follows.
\begin{eqnarray}
w^{N,k}({\cal O}_{h^a}{\cal O}_{h^{b}}|\prod_{i=1}^{n}{\cal O}_{h^{m_i}})_{0,d}&=&
\int_{\widetilde{Mp}_{0,2|n}(N,d)}ev_{0}^{*}(h^a)\cdot ev_{\infty}^{*}(h^b)\cdot\bigr(\prod_{i=1}^{n}ev_{i}^{*}(h^{m_i})\bigr)\cdot c_{top}({\cal E}^{N,k}),\no\\
w^{N,k}({\cal O}_{h^a}|\prod_{i=1}^{n}{\cal O}_{h^{m_i}})_{disk,2d-1}&=&
\int_{\widetilde{Mp}_{D,1|n}(CP^2/RP^2,2d-1)}ev_{0}^{*}(h^a)\cdot\bigr(\prod_{i=1}^{n}ev_{i}^{*}(h^{m_i})\bigr)\cdot c_{top}({\cal E}^{N,k}_{disk}).
\label{mnkdef}
\end{eqnarray}
${\cal E}^{N,k}$ and ${\cal E}_{disk}^{N,k}$ are orbi-bundles. The zero locus of sections of these bundles correspond to quasi maps whose images lie 
inside $M_{N}^{k}$.     
Our discussions in Section 2 and Section 3 are also applicable to these intersection numbers. By combining them with 
the results in \cite{mmg} and \cite{opv}, we obtain the following closed formulas.   
\begin{eqnarray}
w^{N,k}({\cal O}_{h^a}{\cal O}_{h^{b}}|\prod_{i=1}^{n}({\cal O}_{h^{m_i}}))_{0,d}&=&
\frac{1}{(2\pi\sqrt{-1})^{d+1}}\oint_{E^0_{(0)}}\frac{dz_0}{(z_0)^N}\oint_{E^1_{(0)}}\frac{dz_1}{(z_1)^N}\cdots \oint_{E^d_{(0)}}\frac{dz_d}{(z_d)^N}
\times\no\\
&&(z_{0})^a\cdot\biggl(\prod_{j=1}^{d-1}\frac{1}{kz_j(2z_j-z_{j-1}-z_{j+1})}\biggr)\cdot
\biggl(\prod_{j=1}^{d}e^{k}(z_{j-1},z_j)\biggr)\cdot(z_{d})^{b}\times\no\\
&&\prod_{i=1}^{n}\biggl(\sum_{j=1}^{d}w_{m_i}^1(z_{j-1},z_j)\biggr), \;\;(d>0).
\label{mnkth1}
\end{eqnarray}
\begin{eqnarray}
w^{N,k}({\cal O}_{h^a}|\prod_{i=1}^{n}{\cal O}_{h^{m_i}})_{disk,2d-1}&=&
\frac{1}{(2\pi\sqrt{-1})^{d}}\oint_{E^0_{(0)}}\frac{dz_0}{(z_0)^N}\oint_{E^{1}_{(0)}}\frac{dz_1}{(z_1)^N}\cdots \oint_{E^{d-1}_{(0)}}\frac{dz_{d-1}}{(z_{d-1})^N}
2(z_{0})^{a}(k!!)(z_{d-1})^{\frac{k+1}{2}}\times\no\\
&&\biggl(\prod_{j=1}^{d-1}\frac{e^{k}(z_{j-1},z_j)}{(2z_j-z_{j-1}-z_{j+1})}\biggr)\cdot\prod_{i=1}^{n}\biggl(\sum_{j=1}^{d-1}w_{m_i}^1(z_{j-1},z_j)+\frac{1}{2}(z_{d-1})^{m_i-1}\biggr),\no\\
&&\hspace{5cm}(k: \mbox{odd},\;\;d\geq1,\;z_{d}:=-z_{d-1}).
\label{mnkthm2}
\end{eqnarray}
In these formulas, $e^{k}(z,w)$ is the polynomial in $z$ and $w$ given in (\ref{ek}).
These intersection numbers are useful especially in the case of general type hypersurface $M_{N}^{k}$ with $k>N$.
In our previous works, we used the virtual structure constants with two marked points. When the hypersurface is general type, 
we have to operate generalized mirror transformation to translate the virtual structure constants into Gromov-Witten invariants.
According to \cite{gene} and \cite{opv}, it is given as follows.
\begin{eqnarray}
&&w^{N,k}({\cal O}_{h^a}{\cal O}_{h^b})_{0,d}=\no\\
&&=\langle{\cal O}_{h^{a}}{\cal O}_{h^b}\rangle_{0,d}
+\sum_{f=1}^{d-1}\sum_{\sigma_f\in P_f}S(\sigma_{f})
\langle{\cal O}_{h^{a}}{\cal O}_{h^b}\prod_{j=1}^{l(\sigma_{f})}{\cal O}_{h^{1+(k-N)f_j}}\rangle_{0,d-f}
\prod_{j=1}^{l(\sigma_{f})}\frac{w({\cal O}_{h^{N-3-(k-N)f_j}}{\cal O}_{h^0})_{0,f_j}}{k}.\no\\
\label{gmtc}
\end{eqnarray}
\begin{eqnarray}
&&w^{N,k}({\cal O}_{h^a})_{disk,2d-1}=\no\\
&&=\langle{\cal O}_{h^{a}}\rangle_{disk,2d-1}
+\sum_{f=1}^{d-1}\sum_{\sigma_f\in P_f}S(\sigma_{f})
\langle{\cal O}_{h^{a}}\prod_{j=1}^{l(\sigma_{f})}{\cal O}_{h^{1+(k-N)f_j}}\rangle_{disk,2d-2f-1}
\prod_{j=1}^{l(\sigma_{f})}\frac{w({\cal O}_{h^{N-3-(k-N)f_j}}{\cal O}_{h^0})_{0,f_j}}{k}.\no\\
\label{gmto}
\end{eqnarray}
In the above formulas, $P_d$ is a set of usual partition of a positive integer $d$,
\begin{equation}
P_d:=\{\sigma_d=(d_1,d_2,\cdots,d_l)\;|\;d=d_1+d_2+\cdots+d_l,\;1\leq d_1\leq d_2 \leq d_3\leq \cdots\leq d_l\},
\end{equation}
and  $S(\sigma_d)$ is the symmetric factor:
\begin{eqnarray}
S(\sigma_{d}):=\prod_{i=1}^{d}\frac{1}{(\mbox{mul}(i,\sigma_d))!},
\label{symmetric}
\end{eqnarray}   
where $\mbox{mul}(i,\sigma_d)$ is multiplicity of $i$ in $\sigma_{d}$.
Therefore, if we intend to compute $\langle{\cal O}_{h^{a}}{\cal O}_{h^b}\rangle_{0,d}$ and $\langle{\cal O}_{h^{a}}\rangle_{disk,2d-1}$
by using (\ref{gmtc}) and (\ref{gmto}), we have to know the information of the multi point Gromov-Witten invariants $\langle{\cal O}_{h^{a}}{\cal O}_{h^b}\prod_{j=1}^{l(\sigma_{f})}{\cal O}_{h^{1+(k-N)f_j}}\rangle_{0,d-f}$ and $\langle{\cal O}_{h^{a}}\prod_{j=1}^{l(\sigma_{f})}{\cal O}_{h^{1+(k-N)f_j}}\rangle_{disk,2d-2f-1}$ in advance. In the closed string case, we used the associativity equation to compute them. But 
this process made the computation very complicated \cite{gene}. In the open string case, we did not know the open version of the associativity equation 
and we could not compute the open Gromov-Witten invariants genuinely from the open virtual structure constants in the $k>N$ case \cite{opv}.
In contrast, our multi point virtual structure constants should include all the informations of the multi point Gromov-Witten invariants. 
Therefore, we can apply the formalism of Section 2 and Section 3 to execute the generalized mirror transformation for general 
type hypersurface $M_{N}^k$. Let us illustrate our idea from the closed string case. For the hypersurface $M_{N}^{k}$ $k>N$, we introduce 
variables $x^j$ $(j=0,1,\cdots,N-2)$ associated with insertions of ${\cal O}_{h^j}$ and generating functions,
\begin{eqnarray}
&&w^{N,k}({\cal O}_{h^a}{\cal O}_{h^{b}}|(x^0,\cdots,x^{N-2}))_{0}:=\no\\
&&x^c\cdot\int_{CP^{N-1}}k\cdot h^{a+b+c+1}+
\sum_{d>0,m_j\geq 0}w^{N,k}({\cal O}_{h^a}{\cal O}_{h^{b}}|\prod_{j=0}^{N-2}({\cal O}_{h^{j}})^{m_j})_{0,d}\cdot\prod_{j=0}^{N-2}\frac{(x^{j})^{m_j}}{m_j!}.
\label{genchyper2}
\end{eqnarray}
Since (\ref{punc}) and (\ref{divisor}) also hold in this case, these generating functions turn out to be polynomials in $e^{x^1}$ and $x^{j}$ $(j=2,3,\cdots,N-2)$. 
Next, we introduce the generalized mirror transformation,
\begin{equation}
t^{j}(x^0,\cdots,x^{N-2}):=\frac{1}{k}w^{N,k}({\cal O}_{h^{N-2-j}}{\cal O}_{1}|(x^0,\cdots,x^{N-2}))_{0}.
\end{equation} 
If we invert the above equality, our conjecture predicts the following equality.
\begin{eqnarray}
\langle{\cal O}_{h^a}{\cal O}_{h^{b}}(t^0,\cdots,t^{N-2}))_{0}&:=&t^c\cdot\int_{CP^{N-1}}k\cdot h^{a+b+c+1}+
\sum_{d>0,m_j\geq 0}\langle{\cal O}_{h^a}{\cal O}_{h^{b}}\prod_{j=0}^{N-2}({\cal O}_{h^{j}})^{m_j}\rangle_{0,d}\cdot\prod_{j=0}^{N-2}\frac{(t^{j})^{m_j}}{m_j!}=\no\\
&&w^{N,k}({\cal O}_{h^a}{\cal O}_{h^{b}}|(x^0(t^0,\cdots,t^{N-2}),\cdots,x^{N-2}(t^0,\cdots,t^{N-2})))_{0}.
\end{eqnarray}
Let us demonstrate this procedure by taking $M_{8}^{9}$ as an example. In this case, the mirror map is explicitly given up to $d=3$ as follows.
\begin{eqnarray}
t^{0}&=&x^0,\no\\
t^{1}&=&x^1,\no\\
t^{2} &=& x^{2}+34138908q,\no\\
t^{3} &=& x^{3}+124995960x^{2}q+ 8404934443598718q^{2},\no\\
t^{4}&=& x^{4}+249752241x^{3}q+ 
{\displaystyle \frac{340609293}{2}}(x^{2})^{2}q + 
{\displaystyle \frac{123644755203321141}{2}}{x^{2}}q^{2}
+\no\\
&&3815933053700462506215462q^{3}+\cdots,\no\\
t^{5}&=& x^{5}+340609293x^{4}q+556222626
x^{3}x^{2}q +{\displaystyle \frac{257278653}{2}}
(x^{2})^{3}q +\no\\
&&113932607554152477x^{3}q^2 +{\displaystyle \frac{321886193235880779}{2}}(x^{2})^{2}q^2 
+33258838601987300311771653x_{2}q^3+\cdots,\no\\ 
t^{6}&=& x^{6}+374748201x^{5}q +681218586
x^{4}x^{2}q+{\displaystyle \frac{805974867}{2}}(x^{3})^{2}q 
+556222626x^{3}(x^{2})^{2}q +{\displaystyle \frac{257278653}{4}}(x^{2})^{4}q+\no\\
&&139268745219642741x^{4}q^2 
+472782967773195564x^{3}x^{2}q^2 + 
223674801935251734(x^{2})^{3}q^2 
+\no\\
&&48918351923402413916303613x^{3}q^3 
+106098778427559977884727547(x^{2})^{2}q^3+\cdots ,\;\;(q=e^{x^1}).
\label{genemir}
\end{eqnarray}
We can omit the $d\geq 4$ part because Gromov-Witten invarinat of $M_{8}^{9}$ is trivial if $d\geq 4$.
We then compute,
\begin{eqnarray}
\frac{1}{9}w^{8,9}({\cal O}_{h}{\cal O}_{h}|(x^0,\cdots,x^{N-2}))_{0}&=&{x^{4}} + 306470385{x^{3}}q + 215613333(x^{2})^{2}
 + 89761934928094677{x^{2}}q^{2}+ \no\\
&&6297488499797163519141951q^{3}+\cdots. 
\label{whh}
\end{eqnarray} 
If we substitute the inversion of the mirror map to (\ref{whh}), we obtain,
\begin{eqnarray}
\frac{1}{9}\langle{\cal O}_{h}{\cal O}_{h}(t^0,\cdots,t^{N-2})\rangle_{0}&=&t^{4} + ({\displaystyle \frac{90617373}{2}}(t^{2})^{2} 
+ 56718144{t^{3}})e^{t^1} + {\displaystyle \frac{35512880615374365}{2}}{t^{2}}e^{2t^1}+  \no\\
&&1345851991844128981741851e^{3t^1}.
\label{ghh}
\end{eqnarray}
This generating function includes all the non-trivial genus $0$ Gromov-Witten invariants of $M_{8}^{9}$ and reproduce the results in \cite{gene0}!  

In the open string case, we introduce 
\begin{eqnarray}
&&w^{N,k}({\cal O}_{h^a}|(x^0,\cdots,x^{N-2}))_{disk}:=\no\\
&&\sum_{d>0,m_j\geq 0}w^{N,k}({\cal O}_{h^a}|\prod_{j=0}^{N-2}({\cal O}_{h^{j}})^{m_j})_{disk,2d-1}\cdot\prod_{j=0}^{N-2}\frac{(x^{j})^{m_j}}{m_j!}=\no\\
&&\sum_{d>0,m_j\geq 0,\;(j\neq 1) }w^{N,k}({\cal O}_{h^a}|\prod_{j\neq 1}({\cal O}_{h^{j}})^{m_j})_{disk,2d-1}e^{(d-\frac{1}{2})x^1}\cdot\prod_{j\neq 1}
\frac{(x^{j})^{m_j}}{m_j!}.
\label{genchyper}
\end{eqnarray}
Unlike the $CP^2$ case, substitution of inversion of the mirror map to $w^{8,9}({\cal O}_{h^a}|(x^0,\cdots,x^{N-2}))_{disk}$
results in $\langle{\cal O}_{h^a}(t^0,\cdots,t^{N-2})\rangle_{disk}$. 
For examples, we compute,
\begin{eqnarray}
w^{8,9}({\cal O}_{h}|(x^0,\cdots,x^{N-2}))_{disk}&=& (945{x^{2}} + {\displaystyle \frac{945}{2}}{x^{0}}{x^{3}} + 
{\displaystyle \frac{945}{8}}(x^{2})^{2}x^{0}+\cdots)q^{(1/2)} + \no\\
&&(90642729450 +{\displaystyle \frac{236172454245}{2}}{x^{0}}{x^{2}}+\cdots)q^{(3/2)}
+ \no\\
&&({\displaystyle \frac {50109447061228637817}{5}}{x^{0}}+\cdots)q^{(5/2)}+\cdots,
\end{eqnarray}
and
\begin{eqnarray}
w^{8,9}({\cal O}_{1}|(x^0,\cdots,x^{N-2}))_{disk}&=& (945{x^{3}}+{\displaystyle\frac {945}{4}}(x^{2})^{2}+\cdots)q^{(1/2)}+(168225362235x^{2}+\cdots)q^{(3/2)} 
+ \no\\
&&({\displaystyle \frac{276177175032776063634}{25}}+\cdots)q^{(5/2)}+\cdots.
\end{eqnarray}
In this case, we also have the subtlety of infinite insertions of ${\cal O}_{1}$. In the above formula, we truncated the number of ${\cal O}_1$ insertions 
so that the non-trivial open Gromov-Witten invariants of $M_{8}^{9}$ computed in \cite{opv} are not affected.
Substitution of inversion of the mirror map results in,
\begin{eqnarray}
\langle{\cal O}_{h}(t^0,\cdots,t^{N-2})\rangle_{disk}&=&(945{t^{2}}+{\displaystyle \frac{945}{8}}(t^{2})^{2}{t^{0}} + {\displaystyle \frac{945}{2}}{t^{0}}{t^{3}}+\cdots)
e^{(1/2)t^1} + \no\\
&&({\displaystyle \frac{101920638015}{2}}{t^{0}}{t^{2}}+58381461390+\cdots)e^{(3/2)t^1}+ \no\\
&& ({\displaystyle \frac{20865788438073398442}{5}}{t^{0}}+\cdots)e^{(5/2)t^1}, 
\end{eqnarray}
and 
\begin{eqnarray}
\langle{\cal O}_{1}(t^0,\cdots,t^{N-2})\rangle_{disk}&=&(945{t^{3}} + {\displaystyle \frac{945}{4}}(t^{2})^{2}+\cdots)e^{(1/2)t^1} + (33973546005{t^{2}}+\cdots)e^{(3/2)t^1} +\no\\ 
&&({\displaystyle \frac{41731576876146796884}{25}}+\cdots)e^{(5/2)t^1}.
\end{eqnarray}
The integrable condition $\frac{\d}{\d t^0}\langle{\cal O}_{h}(t^0,\cdots,t^{N-2})\rangle_{disk}=\frac{\d}{\d t^1}\langle{\cal O}_{1}(t^0,\cdots,t^{N-2})\rangle_{disk}$
is satisfied. The open Gromov-Witten invariants in these generating functions agree with the results computed from (\ref{ogres}).
\begin{eqnarray}
&&\langle {\cal O}_{h}{\cal O}_{h^{2}} \rangle_{disk,1}=945,\;\langle{\cal O}_{1}{\cal O}_{h^{3}}\rangle_{disk,1}=945,\;
\langle{\cal O}_{1}({\cal O}_{h^{2}})^2\rangle_{disk,1}=945/2,\;\langle{\cal O}_{1}{\cal O}_{h^{2}}\rangle_{disk,3}=33973546005,
\no\\
&&\langle{\cal O}_{h}\rangle_{disk,3}=58381461390,\;
\langle{\cal O}_{1}\rangle_{disk,5}=41731576876146796884/25.
\label{op89data}
\end{eqnarray}
Therefore, we can compute the open Gromov-Witten invariants of $M_{8}^{9}$ genuinely from the open multi-point virtual structure constants. 
Comparing these results with the ones of $CP^2$, the reason why we had to introduce the variables $x^j$ $(j\geq 3) $ in the $CP^2$ case is still unclear. 
We end this paper, leaving pursuit of this subject to future works.

\end{document}